\magnification 1200
\def \sn{{\smallskip \noindent}}
\def \bn {{\bigskip \noindent}}
\def \La {{\longrightarrow}}
\def \la {{\longrightarrow}}
\def \O{{\cal O}}
\def \cO{{\cal O}}
\def \cI{{\cal I}}

\def \Q{{\bf Q}}
\def \bP{{\bf P}}

\def \C{{\cal C}}

\def \F{{\cal F}}

\def \A{{\cal A}}

\def \L{{\cal L}}
\def \o{{\omega}}
\def\bC{{\bf C}}
\font\medium=cmbx10 scaled \magstep1
\font\large=cmbx10 scaled \magstep2

\topskip = 10 truemm

\overfullrule=0pt

\def \C{{\bf C}}

\def \F{{\cal F}}

\def \Q{{\bf Q}}

\def \dim{\mathop {\rm dim}}

\def \noi {\noindent}
\def \dis {\displaystyle}

\def \nmid {\not\>|\>}

\def \cqfd{\unskip\kern 6pt\penalty 500
\raise -2pt\hbox{\vrule\vbox to10pt{\hrule width 4pt
\vfill\hrule}\vrule}\par}

\def\adots{\mathinner{\mkern2mu\raise1pt\hbox{.}
\mkern3mu\raise4pt\hbox{.}\mkern1mu\raise7pt\hbox{.}}}

\def\diagram#1{\def\normalbaselines{\baselineskip=0pt
\lineskip=10pt\lineskiplimit=1pt} \matrix{#1}}
\def\build#1_#2^#3{\mathrel{\mathop{\kern 0pt#1}\limits_{#2}^{#3}}}
\def\hfl#1#2{\smash{\mathop{\hbox to
8mm{\rightarrowfill}}\limits^{\scriptstyle#1}_{
\scriptstyle#2}}}
\def\vfl#1#2{\llap{$\scriptstyle #1$}\left\downarrow \vbox to
4mm{}\right.\rlap{$\scriptstyle #2$}}

\catcode`\@=11

\font\tenmsx=msam10
\font\sevenmsx=msam7
\font\fivemsx=msam5
\font\tenmsy=msbm10
\font\sevenmsy=msbm7
\font\fivemsy=msbm5
\newfam\msxfam
\newfam\msyfam
\textfont\msxfam=\tenmsx  \scriptfont\msxfam=\sevenmsx
 \scriptscriptfont\msxfam=\fivemsx
\textfont\msyfam=\tenmsy  \scriptfont\msyfam=\sevenmsy
 \scriptscriptfont\msyfam=\fivemsy

\def\hexnumber@#1{\ifnum#1<10 \number#1\else
 \ifnum#1=10 A\else\ifnum#1=11 B\else\ifnum#1=12 C\else
 \ifnum#1=13 D\else\ifnum#1=14 E\else\ifnum#1=15 F\fi\fi\fi\fi\fi\fi\fi}

\def\msx@{\hexnumber@\msxfam}
\def\msy@{\hexnumber@\msyfam}
\mathchardef\nmid="3\msy@2D
\mathchardef\varnothing="0\msy@3F
\mathchardef\nexists="0\msy@40
\mathchardef\smallsetminus="2\msy@72
\def\Bbb{\ifmmode\let\next\Bbb@\else
\def\next{\errmessage{Use \string\Bbb\space only in math mode}}\fi\next}
\def\Bbb@#1{{\Bbb@@{#1}}}
\def\Bbb@@#1{\fam\msyfam#1}

\def\bn{\bigskip \noindent}
\def\sn{\smallskip \noindent}
\def\la{\longrightarrow}
\def\cO{{\cal O}}
\def\dim{\rm dim}

\headline={\ifnum\pageno=1 {\hfill} \else{\hss \tenrm -- \folio\ -- \hss} \fi}
\footline={\hfil}

\vskip .5cm 
{\large \centerline {Holomorphic 2-Forms on Complex Threefolds}}

\vskip .5cm \vskip .5cm 
{\medium \centerline {Frederic Campana and Thomas Peternell}}

\vskip .5cm \vskip .5cm

\bn {\medium Introduction}

\bn Holomorphic differential forms play an important role in the global study of  complex projective
or compact K\"ahler manifolds. Forms of degree 1 are well understood: there is a universal map, the Albanese
map, to a complex torus and every 1-form is a pull-back. Forms of top degree also play a special role
because they are sections in a line bundle, the canonical bundle. Forms of degree $d$ with $2 \leq d \leq {\rm dim} X - 1$
are however much harder to study. Especially 2-forms are very interesting, e.g. in symplectic geometry, or by
the fundamental theorem of Kodaira describing them as obstruction for K\"ahler manifolds to be projective. 
In this paper we want to study the "first interesting case" of 2-forms on 3-folds. As general guideline we 
ask whether there is some kind of analogue of the Albanese for 2-forms. This means that we should try to 
find some universal object from which all 2-forms arise. This is impossible to some extend because we can
take products of surfaces with curves and then would have to consider 2-forms on surfaces which are "wild".
So instead we ask whether every 2-forms is produced by some canonical procedure. In sect. 1, which is partially written as
an extended introduction, this is explained in great detail. Cum grano salis, we hope that every 2-form is induced
by a meromorphic map to a surface, to a torus or to a symplectic manifold. We will explain in sect.1, what we 
mean by "induced" (not just pull-back, of course).  
\sn In sections 3,4 and 5 we try to get a picture on 2-forms on 3-folds in terms of the Kodaira dimension. We describe this briefly - a more
detailed account on the results of this paper is given in sect. 2. Let $\omega$ be a 2-form on $X.$
\sn (a) If $\kappa (X) = - \infty $ and if $X$ is not simple (i.e. there is a compact subvariety of positive dimension through 
every point), then there is a meromorphic map onto a surface and $\omega$ is a pull-back.
\sn (b) If $\kappa (X) = 0$ and if $X$ is projective, then after finite cover, \'etale in codimension 2, $X$ is birationally to either
a product of a K3-surface and an elliptic curve or a torus and $\omega$ is induced in an obvious way. The same 
holds also in the K\"ahler case (sect. 8) unless $X$ is simple and not covered by torus, a case which is expected not to
exist.  
\sn (c) If $1 \leq \kappa (X) \leq 2,$ we consider the Iitaka reduction $f: X \la B$ and relate $\omega$ to $f$.
We have to distinguish the cases that $\omega$ is vertical to $f$ or not. In many cases we prove that $\omega$ is
induced.
\sn (d) If $X$ is of general type, we cannot say anything.

\bn We also have more precise results when the zero set of $\omega$ is finite which is in some sense the
"generic" case. 

\bn Sections 6-9 study non-algebraic K\"ahler threefolds $X$. These carry automatically a 2-form. We prove that $X$
is uniruled if $\kappa (X) = - \infty,$ unless $X$ is simple (which case should not exist). 
This is the K\"ahler version of the fundamental result of Miyaoka and Mori in the projective setting.

\bn {\medium 1. Examples and basic notions}
\bn
\rm In this section we give basic examples for projective and K\"ahler threefolds $X$
admitting a non-zero 2-form $\o.$ Our main concern is to try to understand how 2-forms
on a K\"ahler manifold are created, i.e whether the 2-forms come from certain natural
manifolds in analogy to the 1-forms which come from the Albanese torus. The most basic
manifolds admitting 2-forms (in dimension at least 3) are tori and symplectic manifolds, i.e. simply connected manifolds
with an everywhere non-degenerate holomorphic 2-form. 2-forms on surfaces play a completely
different role since they are sections in a line bundle, so their geometry is easy to
understand and there
are plenty of examples.

\bn {\bf 1.1 General Method}
\sn Let $X$ be a compact K\"ahler threefold, $f: X \rightharpoonup S$ a  dominant meromorphic map to  
a surface $S$ with $H^0(K_S) \ne 0$ or a meromorphic map to a torus or a symplectic manifold $S.$
Let $\eta$ be a 2-form
on $S.$ Then $\omega = f^*(\eta)$ is a non-zero 2-form on $X.$ We say that $\omega$ is 
{\it induced} from $S$ or $f.$ More generally we can build $\omega = \sum \omega_i$ where the $\omega_i$
are induced from possibly different maps. 

\bn {\bf 1.2 Vague Question}
\sn {\it Does every 2-form on a K\"ahler threefold arise in this way; i.e. is every 2-form
the sum of induced 2-forms, at least after some finite cover ?}

\bn {\bf Remarks.} (1) One might think of constructing 2-forms as wedge products of 1-forms, say
$\omega = \eta_1 \wedge \eta_2.$ This is in fact a special case of (1.1), i.e. $\omega$ is
induced. Namely, let $\alpha : X \La {\rm Alb}(X)$ be the Albanese map of $X,$ then
$\eta_i = \alpha_i^*(\lambda_i),$ so that $\omega = \alpha^*(\lambda_1 \wedge \lambda_2).$
\sn (2) (1.2) has to be interpreted in the following way in order to avoid obvious counterexamples.
Consider an elliptic fiber bundle $f:X \la S$ over the surface $S.$ Then we can have 2-forms
$\omega$
arising from sections of $\Omega^1_{X\vert S}Ê\otimes f^*(\Omega^1_S).$ If $X = E \times S,$
then $\omega $ is a sum of wedge products of 1-forms on $E$ and $S,$ however if $X$ is
not a product, this is certainly false.

\bn {\bf 1.3 Examples}
\sn (1) Pairs $(X,\omega)$ with $X$ a K\"ahler threefold and $\omega$ a non-zero 2-form exist
in every Kodaira dimension. We just take $X = S \times C,$ the product of a surface $S$ admitting
a 2-form with
a curve $C.$ We can even achieve that $\omega$ has no zeroes unless $X$ is of general type.
For $\kappa \leq 1,$ this
is obvious: just take $S$ to be a torus. For $\kappa = 2,$ take a surface $S$ of general
type admitting a 1-form $\eta$ and a 2-form $\omega'$ without common zeroes. Let $C$ be
elliptic with 1-form $\nu.$ If $p_i$ denote the projections of $X = S \times C,$ we put
$\omega = p_1^*(\omega') + p_1^*(\eta) \wedge p_2^*(\mu).$ Then $\omega$ has no zeroes. 
\sn (2) We now consider fiber products. Let $S_i$ be projective surfaces with maps 
$f_i : S_i \La B$ to the smooth curve $B.$ We assume that the fiber product $S_1 \times _B S_2$
is smooth. This holds if the singular fibers of $f_1$ are over different points than the
singular fibers of $f_2.$ Let $\omega_i$ be 2-forms on $S_i$ and $\alpha_i$ be 1-forms on $S_i.$
Then we obtain a 2-form on $X$ by setting $$\omega = p_1^*(\omega_1) + p_2^*(\omega_2) +
p_1^*(\alpha_1) \wedge p_2^*(\alpha_2).$$
We shall see in (5.5) that every 2-form on $X$ is of this type.
The condition for $\omega$ to be without zeroes is the following:
\sn (a) $\omega_i$ and $\alpha_i$ have no common zeroes, $i = 1,2$
\sn (b) if $f_2^{-1}(b)$ is singular, then $\omega_1$ and $\alpha_1$ have no zeroes along
$f_2^{-1}(b),$ and vice versa for $\omega_2$ and $\alpha_2$ over $f_1^{-1}(b).$ 

\sn (3) Let $S$ be a symplectic 4-fold and $X \subset S$ a smooth threefold. Let $w$ be
the symplectic form on $S$. Then $\omega = w \vert X$ has no zeroes. In fact, choose local
coordinates such that $w = dz_1 \wedge dz_2 + dz_3 \wedge dz_4.$ Now $X$ is locally given 
by $z_4 = f(z_1,z_2,z_3).$ Hence 
$$ w \vert X = dz_1 \wedge dz_2 + dz_3 \wedge {{\partial f}Ê\over {\partial z_1}} dz_1 +
dz_3 \wedge {{\partial f} \over {\partial z_2}} dz_2$$
has no zeroes. 

\sn (4) Let $X$ be a smooth projective threefold and $f: X \La B$ a surjective fibration
to a curve $B.$ Assume that the general fiber of $f$ is an abelian surface and $K_X$ is nef.
Then $X$ carries 2-forms, if $q(X) + h^0(K_X) \geq 2.$ In particular this holds if $g(B) \geq 2.$
In fact, let $F$ be a general smooth fiber of $f.$ First observe       
$$ c_2(X) \cdot F = c_2(F) = 0.$$
On the other hand, $mK_X$ is generated by global sections, $K_X$ being nef. Then either $f$ is the 
Iitaka fibration of $X$ or $mK_X = \O_X.$ In both cases $K_X \equiv \lambda F$ with some
$\lambda \geq 0.$ Thus in total 
$$K_X \cdot c_2(X) = 0,$$
hence $\chi (X,\O_X) = 0$ by Riemann-Roch and the claim on the existence of 2-forms follows. 

\sn (5) (Beauville) Let $S$ be a K3-surface, $S^{(2)}$ its symmetric product and $S^{[2]}$ its
canonical desingularisation, i.e. the blow-up of the image of the diagonal. Then $S^{[2]}$ is
a symplectic 4-fold, call $w$ the symplectic form. If $S$ is projective, Markushevitch [Ma86] proved
the following. If there is a surjective map $\pi: S^{[2]} \La \bP_2$ with general fiber $F$ an
abelian surface, then $S$ admits an elliptic fiber space structure over $\bP_1.$ Consequently
$F$ is a product of elliptic curves. Now let $l \subset \bP_2$ be a line, put $X = \pi^{-1}(l)$ and let $\omega =
w \vert X.$ In this situation $\omega$ is induced. In fact, $X$ is birational to a fiber product
of two elliptic surfaces over $l$ and then $\omega$ is induced as explained in (2).
Therefore one might ask whether all abelian fibrations $X$ over a curve $B$ with $h^2(\O_X) \ne 0$
are birational to a fiber product of two elliptic surfaces. This is however not the case as
explained in the next example.     
\sn (6) (Markushevitch) There is a symplectic 4-fold $M$ (constructed as birational model of
the symmetric product of a K3-surface which is a double covering of $\bP_2$) admitting a surjective map $\pi : M \La
\bP_2$ whose general fiber $F$ is an abelian surface and actually birational to the symmetric 
product $C^{(2)}$ of a hyperelliptic curve of genus 2. In particular $F$ cannot be a product.
As in (5) we let $X = \pi^{-1}(l)$ with $l \subset \bP_2$ a general line, $X =
\pi^{-1}(l)$ and let $\omega = w \vert X.$
Then $\pi: X \La l$ is an abelian fibration over $\bP_1$ with a nowhere vanishing 2-form whose general fiber does
not split. So $X$ is not birational to a product of two elliptic surfaces over $\bP_1.$   

\bn {\bf 1.4 Remarks}  (a) In order to construct further examples where the 2-form is not a priori
induced, one might think of deforming fiber products (1.3(2)) in non-fiber products. 
However it is not possible to get interesting examples (with $\kappa (X) \geq 0,$ as shown in (1.5) below). 
\sn (b) If $f: X \La B$ is an abelian fibration of the projective threefold $X,$ then without
assuming $K_X$ to be nef, we can still produce 2-forms under the condition $q(X) + h^0(K_X) \geq 2.$
Suppose first $g(B) \geq 1.$ Let $X'$ be a minimal model of $X.$ Then $f$ induces an abelian 
fibration $f': X' \La B.$ 
We still have $K_{X'} \cdot c_2(X') = 0.$ However $X'$ might have non-Gorenstein
singularities in which case $\chi(X,\cO_X) > 0,$ the positivity coming from contribution from
the non-Gorenstein singularities, see [Fl87]. Hence $\chi(X,\cO_X) = \chi(X',\cO_{X'}) \geq 0$
and we conclude. If $B = \bP_1,$ we may assume $\kappa (X) \geq 0;$ because if $\kappa (X) = -
\infty,$ we have $q(X) \geq 2$ by our assumption $q(X) + h^0(K_X) \geq 2,$ and therefore we
see immediately that the Albanese map of $X$ is a generic $\bP_1-$fibration over a 2-dimensional
torus, hence we get 2-forms on $X.$ If $\kappa (X) = 0,$ then it is clear by the same 
arguments as sbove that $\chi(X,\cO_X) \geq 0$ and finally if $\kappa (X) = 1, $ we pass to a
minimal model $X'$ and consider its Iitaka fibration $g: X' \La C.$                    

\bn {\bf 1.5 Proposition} {\it Let $g_i : S_i \La C$ be elliptic fibrations of the smooth
compact K\"ahler surfaces $S_i$ onto the smooth curve $B$ such that the fiber product 
$X = S_1 \times _C S_2$ is smooth. Assume $\kappa (S_i) \geq 0.$ Assume furthermore that $g(C) \geq 2$ or 
that the $j-$invariants of $f_i$ are not constant. Then
every small deformation of $X$ is again a fiber product of elliptic surfaces.}

\bn {\bf Proof.} (1) First assume that $g(C) \geq 2.$ Our claim comes down to show that
$$ H^1(T_X) = H^1(T_{S_1}) \oplus H^1(T_{S_2}) \eqno (*) $$
canonically.
Let $p_i : S_1 \times S_2 \La S_i$ be the projections and $f_i = p_i \vert X.$ Let $f: X \La C$
be the canonical map ($f = g_i \circ f_i)$.  
The normal bundle sequence associated to the embedding $X \subset  S_1 \times S_2$ reads
$$ 0 \La T_X \La p_1^*(T_{S_1}) \oplus p_2^*(T_{S_2}) \vert X \La f^*(T_C) \La 0.$$
This sequence induces a map 
$$ \alpha : H^1(T_X) \La H^1(f_1^*(T_{S_1})) \oplus H^1(f_2^*(T_{S_2})) $$
and we want to prove that $\alpha$ is injective with image $H^1(T_{S_1}) \oplus H^1(T_{S_2}).$
The injectivity of $\alpha$ is immediately clear by $H^0(f^*(T_C)) = H^0(T_C)$ by virtue of
our assumption $g(C) \geq 2.$ For the description of the image of $\alpha$ we first compute
by the Leray spectral sequence
$$ H^1(T_X) = H^1(f_{i*}(T_X)) \oplus H^0(R^1f_{i*}(T_X))$$
and
$$ H^1(f_i^*(T_{S_i})) = H^1(T_{S_i}) \oplus H^0(T_{S_i} \otimes R^1f_{i*}(\O_X)).$$
Therefore we must only prove
$$ H^0(T_{S_i} \otimes R^1f_{i*}(\O_X)) = 0. \eqno (**)$$
By base change we have
$$ R^1f_{i*}(\O_X) = g_i^*R^1g_{j*}(\O_{S_j}),$$
where $j=2,$ if $i=1$ and $j=1$ if $i=2.$ Now $\L_j = R^1g_{j*}(\O_{S_j})$ is a seminegative line
bundle, i.e. its dual is nef. In order to prove (**), take a general hyperplane section $B \subset S_i.$
By [Mi87], the vector bundle $\Omega_{S_i}^1 \vert B$ is nef, since $\kappa (S_i) \geq 0.$ So (**) follows
unless possibly ${\rm deg \L_j}Ê= 0$ and if $S_i$ is a torus. This means however that $S_j$ is zero j-invariant. Hence by
our assumption $g(C) \geq 2,$ which contradicts the fact that $S_i$ maps onto $C.$ 

\sn (2) If $g(C) \leq 1,$ we make the same arguments with the following modification. The map $\alpha $ is no longer
injective, however its kernel comes from $H^0(T_C),$ and these elements in $H^1(T_X)$ just correspond to deformations 
of the morphisms $f_i : S_i \La C$ by automorphisms of $C,$ so they also correspond to deformations of $X$ as 
fibered product.

\bn {\it Remark.} Of course (1.5) is false in the case of tori: take $S_i = F_i \times C$ with elliptic curves
$F_i,C$ so that $X = F_1 \times F_2 \times C.$ Then the general deformation of $X$ is no longer a fiber product. 
We shall ignore the case $\kappa (S_i) = - \infty$ in the situation (1.5).

\bn {\bf 1.6 Problem} Let $f: X \La B$ be an abelian fibration of the projective 3-fold $X$ to the smooth  curve
$B.$ Assume possibly $K_X$ nef. Is it true that $c_3(X) \geq 0?$ 
It seems possible to construct examples with $c_3(X) = 0$ by choosing $B$ general in a suitable moduli space,
say the moduli space $\A(1,p)$ of $(1,p)-$polarised abelian varieties with level structures. Note that $c_3(X) = 0$
holds if all fibers of $f$ are of the type $C \times E,$ with $C$ a singular rational fiber of an
elliptic surface and $E$ an elliptic curve.

\bn In the rest of this section we prepare the study of 2-forms on fibered 
threefolds.

\bn {\bf 1.7 Identification}  Let $X$ be a smooth threefold. Then there is a canonical
isomorphism $$\Omega^2_X \simeq T_X \otimes K_X,$$
where $T_X$ is the tangent bundle of $X$ and $K_X$ its canonical bundle. Therefore 
any 2-form $\omega$ corresponds to a section $s \in H^0(T_X \otimes K_X)$ and vice versa and we will
often switch from $\omega$ to $s$ and back.

\bn In the following we fix a compact threefold $X$ and a surjective map $f: X \La B$ to a
smooth curve or a normal surface $B.$   
 
\bn {\bf 1.8 Notation}  Let $Tf: T_X \La f^*(T_B)$ denote the differential with kernel
$T_{X \vert B}.$  Let $s \in H^0(T_X 
\otimes K_X).$  Then we define
$$f_*(s) = f_*(Tf \otimes {\rm id}) (s) \in H^0(B,f_*(f^*(T_B) \otimes K_X)).$$
If $B_0$ is the smooth part of $B,$ then notice that $f_*(s) \vert B_0 \in H^0(B_0,T_{B_0}Ê
\otimes f_*(K_X)).$ 

\bn If $b \in B$ and if $F_b$ denotes the analytic fiber over $b$, then $f_*(f^*(T_B) \otimes
K_X \vert \{b\} = f^*(T_B) \otimes K_X \vert F_b$ and therefore we have
$$ s \vert F_b \in H^0(f^*(T_B) \otimes K_X \vert F_b).$$
By virtue of the exact sequence
$$ 0 \La T_{F_b} \La T_X \vert F_b \La N_{F_b} = f^*(T_B) \vert F_b $$
we deduce that if $f_*(s)(b) = 0,$ then $sÊ\vert F_b \in H^0(F,T_F \otimes K_X \vert F_b).$ 
There is also a morphism
$$ \lambda_b: T_X \otimes K_X \vert F_{\rm red} \La N_{{\rm red} FÊ\vert X} \otimes K_X \vert {\rm red}F.$$
for $F = F_b.$ The condition that $\lambda_b(s) = 0$ is slightly stronger than $f_*(s)(b) = 0.$

\bn {\bf 1.9 Definition} Let $s \in H^0(T_X \otimes K_X).$ Let $b \in B.$ Then we say that
$s$ (or $\omega$) is {\it vertical} at $b,$ if $f_*(s)(b) = 0.$ $s$ is {\it vertical} if it is vertical
at every point.  

\bn In other words, $s$ is vertical if and only if $f_*(s) = 0,$ which is the same as to say
that $s \in H^0(X,T_{X \vert B}Ê\otimes K_X).$  Or, that for x general in $X,$ the annihilator of $\omega$
in $T_{X,x}$ is tangential at $x$ to the fiber of $f$ at $x.$ 

\bn We make the following {\bf observation}.
\sn Suppose $s$ is not vertical at $b.$ Then $f_*(s)$ is a non-zero section of $T_B \otimes f_*(K_X),$ 
which is a torsion free sheaf. Hence $f_*(s) \ne 0$ generically, i.e. $s$ is not vertical at the
general $b \in B.$ 
\sn So $s$ is either vertical everywhere or $s$ is not vertical at the general $b.$ 
\bn We now investigate the following situation. $X$ is a minimal projective 3-fold, i.e. $X$ has
only terminal singularities, is $\Q-$Gorenstein and $K_X$ is nef. Then by abundance $mK_X$ is
generated by global sections for large $m$ and the associated map $f : X \la B$ is the Iitaka
fibration of $X.$ We consider again $s \in H^0(X,T_X \otimes K_X).$

\bn {\bf 1.10 Lemma}  {\it Let $F$ be a fiber of $f.$ Then there exists a neighborhood $V$ of
$f(F)$ in $B$ and a covering $h: \tilde X_V \la X_V = f^{-1}(V),$ \'etale outside a finite set
(contained in ${\rm Sing}(X_V)),$ such that $K_{\tilde X_V}Ê= \cO_{\tilde X_V}$ and such that
$\tilde s = h^*(s) \in H^0(T_{\tilde X_V}).$ }

\bn {\bf Proof.} We have $mK_X = f^*(L).$ Hence $mK_{X_V} = \cO_{X_V},$ if $L \vert V = \cO_V.$
Now take the canonical cover $h: \tilde X_V \la X_V$ with respect to a nowhere vanishing section
of $mK_{X_V},$ see e.g. [KMM87]. Since $K_{\tilde X_V} = \cO_{\tilde X_V},$ we have 
$\tilde s \in H^0(T_{\tilde X_V}).$ 

\bn {\bf 1.11 Notation} In the situation of (1.10) we let $\tilde h : \tilde X_V \la \tilde V$
be the Stein factorisation of $h.$

\bn {\bf 1.12 Corollary}  {\it Assume the situation of (1.10). Let $b \in B$. Then we are in
exactly one of the two following cases.
{\item {(1)} $\omega$ (resp. $s$) is vertical at $b$

\item {(2)} there exists a one parameter group $(g_t)_{t \in \C}$ (after possibly shrinking
$V$) which acts fiber-preserving on $\tilde h : \tilde X_V \la \tilde V$ and which does not 
fix $\tilde b \in \tilde V,$ where $\tilde b $ is any preimage point on $b$ in $\tilde V.$ 
In particular there is a local curve $C_b \subset B $ through $b$ such that $F_{b'} \simeq
 F_b$ for all $b' \in C_b.$}}

\bn {\bf Proof.} Local integration of $s$ near $\tilde h^{-1}(b).$  
\bn \bn 

\bn {\medium 2. Statement of the Main Results}

\bn In this section we give precise statements of the results proved in this paper. 
\bn {\bf 2.0 Notations} By $X$ we denote a compact connected threedimensional KŠhler
manifold. A pair
$(X,\omega)$ denotes a non-zero holomorphic 2-form on the manifold $X$.

\bn {\bf 2.1 Definition} The compact complex threefold $X$ is said
to be {\it
simple} if it contains only finitely many irreducible divisors and only
countably many
compact irreducible curves not contained in any of these divisors (in
particular,
$a(X) = 0$).

\bn In other words, $X$ is simple if there is no positive dimensional compact subspace
through the very general point of $X.$  

\bn {\bf 2.2 Definition} The compact complex threefold $X$ is said
to be {\it  Kummer} if
there exists a surjective meromorphic map $\varphi : T \to X$ with $T$ a torus.

\bn {\bf 2.3 Remark.} The only known examples of simple KŠhler threefolds
$X$ are
Kummer. It can be shown that these should actually be the only ones if a
minimal model
program exists in the KŠhler case when $n=3$ ([Pe96]). Therefore we have the
\bn {\bf 2.4  Conjecture} Every simple KŠhler threefold $X$ is Kummer. 

\bn Since we are far from being able to prove this conjecture, we introduce the

\bn {\bf 2.5 Hypothesis (H)} All compact K\"ahler manifolds are assumed not to be both 
simple and non-Kummer.

\bn {\bf  2.6 Theorem.} {\it Let $X$ be a compact KŠhler threefold with
$\kappa(X) =
-\infty$ $\big($assuming (H)$\big)$. Then  $X$ is uniruled. (see (8.1)).}

\bn This extends the famous projective case ([Mi88],[Mo88]) to the KŠhler case.

\bn {\bf 2.7  Corollary} {\it Let $(X,\omega)$ be as in (2.0) with
$\kappa(X) = -\infty$.
There exists a surjective connected meromorphic map $\rho : X \to S$ to a smooth
surface such that $h^{0,2}(S) > 0$ and $\rho^\ast : H^{0,2}(S) \to
H^{0,2}(X)$ is bijective.}

\bn One just needs to take the rational quotient ([Ca81,92]) of $X$ to get $\rho$, see
(3.1) for details. Recall also that non-algebraic K\"ahler threefolds carry a holomorphic
2-form by Kodaira's theorem. 
A classification of non-algebraic uniruled threefolds $X$, according
to $a(X),$  is given
in sect. 9.

\bn {\bf 2.8 Theorem}  {\it Let $X$ be a compact KŠhler threefold with
$\kappa(X) = 0$
and $h^{0,2}(X) > 0$ (assuming (H)). Then $X$ admits a bimeromorphic model
$X'$ with a
covering $\pi : \widetilde X' \to X'$, Žtale in codimension one, such that
$\widetilde X'$ is
either a torus or a product $S\times E$ with $S$ a K3-surface and $E$ an
elliptic curve
(see 8.1).}

\bn Again, this extends to the KŠhler case a known result ([Ka81]) in the
projective case.

\bn {\bf 2.9 Corollary} {\it In the situation of (2.8) let $\tilde \omega$ be 
the form on $\tilde X',$ induced by $\omega.$  
Then
$\tilde \omega$ is induced from maps to either tori or K3-surfaces.}

\bn When $\kappa \ge 1$, the situation becomes more complicated, and our
results are
very partial. When $\kappa(X) = 3$, we shall not say anything. On the other
hand, we
don't need any longer the hypothesis (H), since $a(X) \ge \kappa(X) \ge 1$.
We shall
denote in 2.10 and 2.11 below by $f : X \to B$ the (meromorphic) Iitaka reduction of $X$.

\bn {\bf 2.10 Theorem} {\it Let $(X,\omega)$ be as in 2.0 with 
$\kappa(X) =
1$. Then:\hfill\break
{\rm (i)} If $\omega$ is not $f$-vertical, then $f$ is generically over $B$
a meromorphic
fiber bundle (i.e. two generic fibers are bimeromorphic). When $X$ is
projective, more
precise conclusions can be drawn (see 4.1, 4.2).\hfill\break
{\rm (ii)} If $\omega$ is $f$-vertical, and $X$ is {\it minimal
projective}, the generic fiber
$F$ of $f$ is either abelian or bielliptic. If $F$ is not a simple abelian
surface, then after a
finite base change $\beta : \overline B \to B$, the space $\overline X := X \dis
\mathop
{\times}_{B} \overline B$ is a fiber product $\overline X = X_1 \dis \mathop
{\times}_{\overline B} S_2$ of two elliptic surfaces and $\overline
\omega$ (the
lift of $\omega$ to $\overline X$) is induced. }

\bn  See (4.4), (4.5) for details and further conclusions in even more special cases.

\bn {\bf 2.11 Theorem.} {Let $(X,\omega)$ be as in 2.0, and assume
that $\kappa(X)
= 2$ with $X$ a projective ${\Bbb Q}$-factorial threefold with only
terminal singularities
and $K_X$ nef (so that $f$ is holomorphic with $B$ a normal surface with
only quotient
singularities).\hfill\break
{\rm (i)} If $\omega$ is $f$-vertical, then $\omega = f^\ast(\eta)$ where
$\eta$ is a
$2$-form on $B$ (see 5.1).
\hfill\break
{\rm (ii)} If $\omega$ is not $f$-vertical and the invariant $j : B \to
{\Bbb P}_1$ is not
constant, then $X$ is birational to $S \dis \mathop {\times}_{C}
B$ where $J : B
\to C$ is the Stein factorisation of $j$ and $\varphi : S \to C$ a suitable
elliptic fibration.
Moreover, $\omega$ is induced from $B$ and $S$.}

\noi Let us now briefly explain the proofs of our main results 2.6 and 2.8
: the results are
known when $a(X) = 3$ by the results on the minimal model program (Kawamata,
Miyaoka, Mori). The case $a(X) = 0$ is again very easy (under assumption
(H)). We proceed
by considering the algebraic reduction $r : X \to S$ in the remaining cases
$a(X) = 2, 1$ :
the generic fiber $F$ or $r$ is either an elliptic curve; a ruled surface
(the $X$ is clearly uniruled) or a
surface with $\kappa(F) = 0$, in fact either K3 or torus (bimeromorphically). The
difficult
cases are the cases where $F$ is a torus (either an elliptic curve or
2-dimensional).

\noi Our main tools are first the solutions of Iitaka conjectures
$C_{3,1}$, $C_{3,2}$
([Fu78], [Ka81]) and secondly the {\it ``untwisted model"} (see ¤ 6) $r_0 : X_0 \to S$ of $r$, whose
generic fiber is Aut$^0(F)$, the group of translations of $F$. The fiber
$F_0$ of $r_0$
corresponding to $F$ is (non-canonically$\,!$) isomorphic to $F$, and $r_0$ has
moreover a meromorphic canonical zero section. When ${\dim} X = 2 = 1 +
{\dim} S$, this is
the classical Jacobian fibration used by K. Kodaira in his classification theory of surfaces.

\noi We show that $h^p({\cal O}_{X_0}) = h^p({\cal O}_X)$ $(p \ge 0)$ and
$\kappa(X_0) \le
\kappa(X)$ (in special cases, sufficient for our purposes, but this should
be a general fact).

\noi To illustrate the method, let us explain how it works when $a(X) = 2$,
$\kappa(X) =
-\infty$.

\noi If $\kappa(S) \ge 0$, we are finished by $C_{3,2}$ which shows that
the generic fiber
of $r$ is ${\Bbb P}_1$.

\noi Otherwise, we consider the untwisted model $r_0: X_0 \to S.$ Since $r_0$ has a section, $X_0$ is projective.
By general theory (6.5) we have $\kappa (X_0) = - \infty,$ therefore $X_0$ is uniruled and we consider its
rational quotient $X_0 \to \Sigma_0$ whose image turns out to be a surface. From this map
we construct another map $\sigma : X \to \Sigma$ with $\Sigma$ a
surface with $\kappa(\Sigma) \ge 0$ (since $h^{0,2}(X) > 0$)
(see ¤ 8).

\noi It might be noted that most (not all) of our arguments have an
inductive nature on
${\dim} X$, and also it seems quite plausible that if the implication
``$\kappa(X) =
-\infty$" implies ``$X$ uniruled" can be proved in the projective case, it
should also hold
in the KŠhler case as well (provided $C_{n,m}$ holds true, too).

\bn \bn \bn   {\medium 3. The Case of non-positive Kodaira Dimension}

\bn For all of this section we fix a compact K\"ahler threefold $X$ and a non-zero 2- form $\omega$ on $X.$ 
We shall investigate the structure of $(X,\omega)$ in case $\kappa (X) \leq 0.$ 

\bn {\bf 3.1 Theorem} {\it If $\kappa (X) = - \infty,$ then there exists a meromorphic dominant map 
$f: X \rightharpoonup S$
to a K\"ahler surface $S$ such that $\omega = f^*(\eta),$ unless $X$ is simple.}

\bn Recall that $X$ is {\rm simple}, if there is no covering family of proper positive-dimensional subvarities, in particular
$X$ has non non-constant meromorphic function. It is expected that there are no simple threefolds of negative Kodaira
dimension, see sect. 2. 

\bn {\bf Proof.} (a) First note that $X$ is uniruled. In the algebraic situation this is a central result in minimal model
 theory; if $X$ is K\"ahler with $a(X) < 3,$  this is Theorem 8.1, since $X$ is supposed not to be simple. 
  
\sn (b) Since $X$ is uniruled, there exists a covering family $(C_t)$ of rational curves and we can consider an associated rational
quotient $f: X \rightharpoonup S,$ i.e. $f$ contracts just the curves $C_t.$ Since the general fiber 
of $f$ is rational, we find a
2-form $\eta$ with $\omega = f^*(\eta).$ In fact, choose a sequence of blow-ups $\hat X \La X$
such that the induced map $\hat f: \hat X \La S$ is a morphism. Let $S_0$ be maximal such that
$\hat f $ is a submersion over $S_0,$ i.e. a $\bP_1-$bundle. Then we find a 2-form $\eta_0$ on
$S_0$ such that $\omega = \hat f^*(\eta_0)$ over $X_0 = \hat f^{-1}(S_0).$ Since $f^*(\eta_0)$
extends as a 2-form on $X,$ it is easily checked that 
$$ \int_{S_0}\eta_0 \wedge {\overline {\eta}}_0 < \infty,$$
hence $\eta_0$ extends to all of $S.$ 

\bn {\bf 3.2 Theorem}  {\it Assume that $X$ is projective. If $\kappa (X) = 0,$ then after possibly a finite cover, \'etale in codimension 2, 
either $X$ is bimeromorphic to a torus, $\omega$ coming from this torus, or 
$X$ is birationally equivalent to $E \times S,$ with $E$ an elliptic curve and $S$ a K3 surface, $\omega$ coming from the K3 surface}

\bn {\bf Proof.} (a) First assume that $q(X) > 0.$ By passing to a minimal
model we may assume from the beginning that $K_X \equiv 0.$ Of course, $X$ is now singular in
general. Let $\alpha : X \La A$ be the Albanese map. After performing
a base change, \'etale in codimension 2, $\tilde A \La A,$ $X$ is by [Ka85] birationally a 
product $F \times A$ with $F$ a K3 surface or empty.      
\sn (b) We now prove that the existence of a 2-form implies the existence of a 1-form if $K_X = \cO_X$ which 
can be achieved after finite cover, \'etale in codimension 2.  
In fact, by Riemann-Roch [Fl87], we have $\chi(X,\cO_X) = 0.$ So $q(X) > 0$ and (a) applies. Compare [Pe94,sect.5] for
an argument without using Riemann-Roch (and which applies directly also in the K\"ahler case).  
    
\bn We are now turning to the special case that the 2-form $\omega$ has only finitely many zeroes. 

\bn {\bf 3.3 Theorem}  {\it Let $X$ be a projective threefold and $\omega$ a 2-form on $X$
with finite zero set $Z(\omega).$ 
{\item {(1)}Assume $\kappa (X) = - \infty.$ Then $X$ is a $\bP_1-$bundle over a K3-surface
or a torus and $\omega$ is a pull-back.
\item {(2)} If $\kappa (X) \geq 0,$ then $K_X$ is nef.} }

\bn {\bf Proof.} Assume $K_X$ not to be nef. If $\kappa (X) = - \infty,$ this is automatically
satisfied.
Let $\varphi : X \La Y$ be the contraction of an extremal ray. Let $d =  {\rm dim}Y.$ 
Since $H^2(X,\O_X) = H^2(Y,\O_Y),$ we must have $d \geq 2.$ 
\sn We first treat the case $d = 2.$ Here $X \La Y$ is a conic bundle over the smooth surface $Y.$
Let $\Delta \subset Y$
be the discriminant locus. If $\Delta = \emptyset,$ then $\varphi$ is analytically a $\bP_1-$
bundle and there is a 2-form $\eta$ on $S$ with $\omega = \varphi ^*(\eta).$ Since $\eta$ has
no zeroes, $S$ is a torus or a K3-surface as claimed. So suppose $\Delta \ne \emptyset.$ 
Pick a smooth point $y \in \Delta.$ Then $\varphi ^{-1}(y)$ is a  reducible conic in $\bP_2,$
i.e. consists of 2 lines $L$ and $L'$ with normal bundles (in $X$) $\O \oplus \O(-1).$ 
By taking $\bigwedge ^2$ of the exact sequence
$$ N^*_{L \vert X} \La \Omega^1_{X} \vert L \La \Omega^1_L \La 0$$
we obtain the exact sequence
$$ 0 \La \bigwedge ^2 N^*_{L \vert X} \La \Omega^2_X \vert L \La \Omega^1_L \otimes N^*_{L \vert
X} \La 0,$$
which reads
$$ 0 \La \O_L(1) \La \Omega^2_X \vert L \La \O_L(-2) \oplus \O_L(-1) \La 0.$$
Consequently $\omega$ vanishes at some point of $L.$ Varying $y$ we conclude ${\rm dim} Z(\omega) \geq 1,$
contradiction. Hence $\Delta = \emptyset.$
\sn (3) We finally consider the case $d = 3,$ so that $\varphi$ is birational. Here we use
Mori's classification [Mo82]. If $\varphi$ is the blow-up of a smooth curve in the smooth
threefold $Y,$ then the abover arguments applied to a non-trivial fiber of $\varphi$ give
the same contradiction. Let $E$ be the exceptional divisor of $\varphi$. Then we conclude
${\rm dim}\varphi (E) = 0.$ If $E = \bP_1 \times \bP_1$ with normal bundle $\O(-1) \oplus \O(-1),$
then we have a factorisation $\varphi = \rho \circ \psi,$ where $\psi : X \La Z$ is the
blow-down of $X$ along one of the two rulings of $E$ and $\rho$ is a small contraction.
Then we apply the above arguments to $\psi$ and exclude that case.

Next suppose that $E = \bP_2$ with normal bundle $\O(-1),$ i.e. $\varphi$ is the blow-up of
a smooth point in $Y.$ We have the exact sequences
$$ 0 \La N^*_{E \vert X}Ê\La \Omega^1_X \vert E \La \Omega^1_F \La 0$$
and 
$$ 0 \La N^*_{E \vert X} \otimes \Omega^1_F \La \Omega^2_F \La 0.$$
Therefore $H^0(\Omega^2_X \vert F) = H^0(N^*_{F \vert X}) = H^0(\Omega^1_F(1)),$ and the
last group vanishes since $ E = \bP_2.$ So this case is not possible, too. 
The same arguments also exclude $E = \bP_2$ with normal bundle $\O(-2).$ 

It only remains to consider the case that $E$ is a singular normal quadric in $\bP_2$ with
normal bundle $N_{E \vert X}Ê= \O(-1).$ Let $x_0$ be the vertex of $E.$
The exact sequence 
$$0 \La N^*_E \La \Omega^1_X \vert E \La \Omega^1_E \La 0$$
gives an exact sequence 
$$ 0 \La N^*_E \otimes \Omega^1_E \buildrel \alpha \over \La \Omega^2_X \vert E  \buildrel
\beta \over \La \Omega^2_E \La 0$$
at least on $E \setminus \{x_0\}.$ 
Here $\Omega^2_E := \bigwedge^2 \Omega^1_E / {\rm torsion}$. Actually it can be seen that
$\bigwedge^2 \Omega_E$ is torsion free, but we do not need that. 
Now the maps $\alpha$ and $\beta$ extend to maps
$$ N^*_E \otimes \Omega^1_E \La \Omega^1_X \vert $$
and $$ \Omega^1_X \vert E \La \Omega^2_E.$$
The second extension being obvious, the first comes from the reflexivity of $\Omega^1_E.$ 
It is now immediately checked that the second sequence remains exact on all of $E.$ 
Since $\Omega^2_E \subset \omega_E,$ the dualising sheaf of $E$ we have $H^0(\Omega^2_E) = 0,$
hence it suffices to prove $$H^0(E,\Omega^1_E(1)) = 0$$
to obtain the desired contradiction. This is however obvious from the embedding $E \subset 
\bP_3.$

\bn In the non-algebraic case we can prove

\bn {\bf 3.4 Theorem}  {\it Let $X$ be a compact K\"ahler threefold with $\kappa (X) = - \infty.$
Assume that $X$ is not simple and that there is a 2-form $\omega$ on $X$ with $Z(\omega)$ finite.
Then $X$ is a $\bP_1-$bundle over a torus or a K3-surface.}

\bn {\bf Proof.} By Theorem 8.1, $X$ is uniruled. Therefore [Pe96] applies and there is a
contraction $\varphi : X \La Y$. Hence we can conclude as in 3.3.

\bn In the proof of 3.3 the assumption $\kappa (X) = - \infty$ was only used to ensure the existence of
a contraction. Therefore we can state more generally :

\bn {\bf 3.5 Corollary}  {\it Let $X$ be a projective threefold with a 2-form having at most
finitely many zeroes. If $K_X$ is not nef, then $X$ is a $\bP_1-$bundle over a torus or a
K3-surface.}

\bn {\bf 3.6 Theorem} {\it Let $X$ be a projective threefold with $\kappa (X) = 0$ and $\omega$
a 2-form with finite zero set. Then $X$ is a torus or a product of a K3-surface and an elliptic
curve up to finite \'etale cover.}

\bn {\bf Proof.} $K_X$ is nef by 3.5. Since $\kappa (X) = 0,$ it follows $K_X \equiv 0.$ 
Then the result follows from the decomposition theorem for manifolds with $c_1(X) = 0,$
observing that Calabi-Yau threefolds have no 2-forms. 

\bn In the non-algebraic case we cannot conclude since we do not know the existence of a contraction when $K_X$ is not nef.
However it should be possible to prove (3.6) for non-simple K\"ahler threefolds with $\kappa (X) = 0.$
By [Fu83], we have $a(X) \geq 1.$ Then one should analyse the algebraic reduction.

\bn \bn 
{\medium 4. Case $\kappa = 1.$}

\bn Let $X$ be a compact K\"ahler 3-fold with 2-form $\omega.$ Throughout this section we
assume $\kappa (X) = 1.$ Let $f: X \rightharpoonup B$ be the Iitaka fibration to the smooth
curve $B.$ By possibly blowing up $X,$ we may assume $f$ holomorphic. Anyway $f$ is automatically
holomorphic if $g(B) \geq 1.$
We consider the exact sequence 

$$ 0 \la T_{X \vert B}Ê\otimes K_X \la T_X \otimes K_X \buildrel \kappa \over \la 
f^*(T_B) \otimes K_X.$$
Then  $\omega $ resp. $s$ is vertical iff $\kappa (\omega) = 0.$ 
Equivalently, consider 
$$ f^*(\Omega^1_B) \otimes \Omega^1_{X \vert B} \la \Omega^2_X \buildrel \mu \over \la \Omega^2_{X \vert B}
\la 0.$$ Then $\omega $ is vertical iff $\mu (\omega) = 0$ generically. So if $\omega$ is not vertical,
then it defines an element $\overline {\omega}Ê\in H^2(X,\cO_X),$ and since $\mu (\omega) 
\ne 0,$ we have $0 \ne \overline {\omega}Ê\vert X_b \in H^2(X_b,\O),$ at least for general $b.$ 
Therefore we obtain a section
$$ t \in H^0(B,R^2f_*(\O_X)).$$ 
Since $f_*(\omega_{X \vert B}) \simeq R^2f_*(\cO_X)^*$ is nef, we conclude
$$ f_*(\omega_{X \vert B})Ê\simeq \O_B,$$
thus by Torelli (we may assume all fibers minimal, see e.g. [Ue87,1.11]), $f: X_U \la U$
is an analytic fiber bundle, where $U$ is the largest open set over which $f$ is smooth.

\bn We conclude

\bn {\bf 4.1 Proposition}  {\it Let $X$ be a compact K\"ahler 3-fold with 2-form $\omega$ and
with holomorphic Iitaka fibration $f: X \la B.$ Assume that $\omega$ is not vertical with
respect to $f.$ Then the smooth part $f: X_UÊ\La U$ is an analytic fiber bundle.}

\bn  For further investigation let $X$ be projective and $X'$ be a minimal model of $X.$
Let $f': X' \la B$ be the Iitaka fibration, i.e. the fiber space defined by $\vert mK_{X'}Ê\vert$
for suitable $m.$ Again we let $s \in H^0(X,T_X \otimes K_X)$ be the section associated to 
$\omega.$  Then the section in 
$H^0(f'_*(\omega_{X' \vert B}))$ associated to $\omega$ has no zeroes as seen above, hence this
sheaf is trivial. Therefore
$s$ does not vanish on full fibers of $f'.$ Now pass to $\tilde X'_V \la \tilde V$ according to
(1.10) and (1.12). Then actually $\tilde s \in H^0(T_{\tilde X'_V}).$ Since $\tilde s$ is
not vertical at every $\tilde b \in \tilde B$ by the above non-vanishing, (1.12) implies that
all fibers of $\tilde f$ are smooth (singular fibers would be move horizontically to other
singular fibers). We conclude that $f': X' \setminus {\rm Sing}(X) \la B$ is almost smooth, i.e.
the only singular fibers are multiples of smooth surfaces. Now let $b_0 \in B$ such that 
$$  F_{b_0}Ê\cap {\rm Sing}(X') \ne \emptyset, $$
where $F_{b_0}Ê= f'^{-1}(b_0).$ Then $S = {\rm red}F_{b_0}$ is a normal surface with a certain
fiber multiplicity $\lambda \geq 1.$ If $\tau : \tilde V \la V$ is the canonical map as in 
(1.12), we take $\tilde b_0 \in \tilde V$ such that $\tau (\tilde b_0) = b_0.$ Then $\tilde
S = {\rm red}F_{\tilde b_0}$ is smooth, and we have a birational finite map $\rho : \tilde S
\la S.$ Since $S$ is normal, $\rho$ is biholomorphic, and $S$ is smooth. Therefore $f'$ is
almost smooth. Moreover there is a finite base change $\tilde B \la B$ such that the induced
fiber space $\tilde f: \tilde X \la \tilde B$ is smooth, hence an analytice fiber bundle.
Note that the typical fiber must be a K3 surface or abelian, since $F$ is minimal with $K_F$
numerically trivial; it cannot be hyperelliptic or an Enriques surface because of our 2-form
$\omega.$ 
If $F$ is a K3 surface, then, ${\rm Aut}(F)$ being discrete, $\tilde X = F \times \tilde B$
after possibly another finite \'etale cover. So assume $F$ abelian. It is well known (see e.g.
[Fu83]) that after passing to a finite \'etale cover of $\tilde B,$ we have
$$ q(\tilde X) = q(F) + q(\tilde B).$$ 
Now consider the Albanese map $\tilde \alpha : \tilde X \la \tilde B.$ After a last \'etale
base change we will have $\tilde X = F \times \tilde B,$ see [Ue75] (note that since $\kappa (\tilde X) \geq 1,$ we
must have $g(\tilde B) \geq 2)$.

\bn In total we have proved

\bn {\bf 4.2 Theorem}  {\it Let $X$ be a smooth projective threefold with 2-form $\omega.$
Assume that $\omega$ is not vertical with respect to the Iitaka fibration. Let $X'$ be a
minimal model with Iitaka fibration $f': X' \la B.$ Then $f'$ is almost smooth. There is a
finite base change $\tilde B \la B$ with induced fiber space $\tilde f : \tilde X \la \tilde B$
such that $\tilde X \simeq F \times \tilde B,$ where $F$ is abelian or K3. If $\omega'$ is the
induced 2-form on $X'$ and $\tilde \omega $ the pull-back to $\tilde X,$ then $\tilde \omega =
p^*(\eta) + p^*(\eta) \wedge \tilde f^*(u),$ where $p: \tilde X \la F$ is the projection and $u,v$
are 1-forms.}

\bn Of course (4.2) holds in the K\"ahler case as well since also in this case the existence of a minimal
model is known (Nakayama [Na95]).

\bn {\bf (4.3)} In case that $\omega$ is vertical with respect to the Iitaka fibration we cannot 
say so much. The example 1.3(6) shows that the Iitaka map is not almost smooth in general for
smooth minimal threefolds. We make this more precise and take over the notations of 1.3(6).
We let $f = \piÊ\vert X.$ Then $K_X = f^*(\cO_l(2)),$ hence $\kappa (X) = 1$, $K_X$ is nef and
$f$ is the Iitaka fibration. Assume $f$ almost smooth. Then $$f_*(\omega_{X \vert l}) = \cO_l,$$
so dually $R^2f_*(\cO_X) = \cO_l.$ From the Leray spectral sequence we deduce $h^2(\cO_X) \geq 2.$ 
On the other hand we have the exact sequence
$$ \bC \simeq H^2(\cO_M) \la H^2(\cO_X) \la H^3(\cI_X) \la H^3(\cO_M) = 0.$$
Now $H^3(\cI_X) = H^1(\cO_M(X)) = H^1(\pi^*(\cO_{\bP_2}(1))).$ 
By the Leray spectral sequence
$$ H^1(\pi^*(\cO_{\bP_2}(1))) \subset H^1(\cO_{\bP_2}(1)) \oplus H^0(R^1\pi_*(\cO_M) \otimes
\cO_{\bP_2}(1)).$$
Since $q(M) = 0,$ the Leray spectral sequence immediately implies that $R^1\pi_*(\cO_M)$ is a
negative vector bundle, hence $H^3(\cI_X) = 0$ and $h^2(\cO_X) = 1,$ contradiction.
So $f$ is not almost smooth. 

\bn {\bf 4.4 Proposition} {\it Let $X$ be a minimal K\"ahler threefold with $\kappa (X) = 1.$ Let
$\omega$ be a 2-form on $X$ which is vertical with respect to the Iitaka fibration
$f: X \la B.$ Let $F$ be the general fiber of $f.$
{\item {(1)} $F$ is birationally an abelian surface or hyperelliptic.
\item {(2)} If $F$ is hyperelliptic or a non-simple abelian surface, then there exists a finite cover $\tilde X \la X,$ 
unramified over the largest open set $B_0$ over which $f$ is smooth, such that $ \tilde f: \tilde X \la \tilde B$
is birational to a fiber product $S_1 \times _{\tilde B} S_2.$ 
\item{(3)} Assume $q(X) > q(B).$ Then $F$ is hyperelliptic or a non-simple abelian surface, hence 
(2) holds, or $X$ is birational to an torus fiber bundle over a curve of genus $g \geq 2$ which gets a product
after finite \'etale cover.
\item{(4)} If $B$ is rational, then $q(X) = 1$ and the general fiber $G$ of the Albanese map has $\kappa (G) = 1.$} }

\bn {\bf Proof.}  (1) Since $\omega$ is vertical, it defines a 1-form on every smooth fiber $F$ of $f.$  Since
$\kappa (F) = 0,$ $F$ must be abelian or hyperelliptic.
\sn (2) If $F$ is hyperelliptic, then $F$ becomes a product after finite \'etale cover. The same holds if $F$ is
a non-simple abelian surface by Poincar\'e complete reducibility, since then $F$ has the structure of an elliptic
fiber bundle over an elliptic curve. Now consider the family of elliptic curves from
the projection maps of the smooth fibers $F$ and take closure in the Chow scheme. After a finite
cover of the base, unramified over $B_0,$ we obtain a second family of elliptic curves and
therefore the cover $\tilde X \la X$ is birationally a fiber product over $\tilde B.$ The technical
details are left to the reader. 
\sn (3) Suppose now that $q(X) > q(B).$ We have only to show that if $F$ is abelian, then either it
is not simple or we are birationally in a product situation after some \'etale cover. 
Let $\alpha : X \la {\rm Alb}(X)$ be the Albanese map of $X$ and $\beta: 
B \la {\rm Alb}(B)$ that one of $B.$ We obtain a commutative diagram
$$  \matrix{ X & \la & {\rm Alb}(X) \cr
           \downarrow &  & \downarrow \cr
            B & \la & {\rm Alb}(B) \cr }.$$
If ${\rm dim} \alpha (F) = 0,$ then $\alpha$ would factor over $f$ which contradicts
$q(X) > q(B).$ If ${\rm dim}Ê\alpha (F) = 1,$ Then $F$ is not simple. So suppose ${\rm dim} \alpha (F) = 2.$ Then $\alpha (F)$ is a subtorus of
${\rm Alb}(X).$ If $\alpha (X) = \alpha (F),$ we have $q(X) = 2,$ both 2-forms coming from
$F.$ Hence $q(B) = 0.$ This case is excluded in (4). So we are reduced to the case
${\rm dim} \alpha (X) = 3.$ Now $\alpha \vert F$ is \'etale, hence an isomorphism after an
\'etale base change. This implies by the diagram that $\alpha$ is birational. Therefore
we may substitute for our purpose $X$ by $Y$.
Since $\kappa (X) = 1,$ the map $Y \la \beta (B)$ is an \'etale fiber bundle [Ue75] and hence
a product after finite \'etale cover. 

\sn (4)  Assume that $B$ is rational. Let $\alpha : X \la Y \subset A$ denote again the Albanese map
of $X.$ Let $G$ be a general fiber of $\alpha.$ Suppose $\dim Y = 1.$ Then $C_{3,1}$ yields $\kappa (G)
= 0$ unless $Y = A$ is an elliptic curve. If $\kappa (G) = 0,$ then by nefness of $K_X$, we deduce
$K_G \equiv 0$ and $K_X \cdot G = 0.$ Since $f$ is the Iitaka fibration, $\dim f(G) = 0,$ hence there
is a factorisation $X \la Y \la B$. Therefore $f$ does not have connected fibers, contradiction. 
So $Y$ is elliptic and $\kappa (G) = 1$ as claimed.
\sn Thus we may assume $\dim Y \geq 2.$ If $\dim \alpha (F) = 2,$ then we obtain a 2-form on $X$ which
is not vertical with respect to $f,$ hence $f$ is almost holomorphic by (4.2). Thus $\kappa (X)  \ne 1$
(or argue as follows: the sheaves $R^{i}f_*(\cO_X)$ are trivial, hence the Leray spectral sequence 
gives $q(X) = 2$ which excludes all cases except $Y = {\rm Alb}(X)$). It remains to treat the case
that $\dim \alpha (F) = 1.$ If $Y = {\rm Alb}(X),$ then $A$ is not simple, we have a bundle structure
$A \la C$ to an elliptic curve, contracting all the curves $\alpha (F),$ and $f$ factors as $X \la Y
\la C \la B.$ Therefore the fiber of $f$ cannot be connected, contradiction. If $Y \ne A,$ then 
there is a map $Y \la C$ to a curve of general type, contracting again all $\alpha (F),$ and we 
argue as before.

\bn We finally have a short look at the case that $\omega$ has only finitely many zeroes.

\bn {\bf 4.5 Proposition} {\it Let $X$ be a smooth K\"ahler threefold with $\kappa (X) = 1.$
Let $\omega $ be a 2-form with only finitely many zeroes, not vertical with respect to the
Iitaka fibration. Then $K_X$ is nef, therefore we
have a holomophic Iitaka fibration $f: X \la B$. $f$ is almost holomorphic. There is a finite cover $\tilde B \la B$
such that the induced fiber space $\tilde f: \tilde X \la \tilde B$ splits:
$\tilde X \simeq \tilde B \times S$ with $S$ a K3 surface or abelian, $\tilde \omega$ coming from
$S.$ In particular $\omega$ has no zeroes.} 

\bn {\bf Proof.}  $K_X$ is nef by (3.3). Therefore $mK_X$ is generated by global sections for
suitable $m,$ and the Iitaka fibration is holomorphic. The rest comes from the proof of (4.2). 

\bn \bn 
{\medium 5. Case $\kappa = 2.$} 

\bn We now turn to the case that $\kappa (X) = 2.$ Throughout this section we assume that
$X$ is a K\"ahler ${\bf Q-}$factorial threefold with only terminal singularities such that
$K_X$ is nef; so $X$ is a minimal model. Then the Iitaka fibration is automatically a
holomorphic elliptic fibration, and the surface $B$ has at most quotient singularities. 
We let $S \subset B$ be the closure of the set of all $b \in B$ such that the fiber
$f^{-1}(b)$ is a multiple elliptic curve. Let $S_i$ be the irreducible components of $S_i.$
\sn For any space $Z$ we let $\Omega^{i}_Z = \bigwedge_{i} \Omega^1_Z$ and if $Z$ is normal
we define $\tilde \Omega^{i}_Z = j_*(\Omega^{i}_{{\rm reg} Z}).$ 

\bn {\bf 5.1 Lemma}  {\it Assume that the 2-form $\omega$ is vertical with respect to $f.$ 
Then $\omega = f^*(\eta)$ with $\eta \in  \omega_B(\sum k_iS_i) = \Omega^2_B(\sum k_iS_i)$ 
with $k_i \geq 0.$}

\sn {\bf Proof.} 
Let $B_0$ be the largest subset of $B$ such that $f$ is a submersion over $B_0.$ Let
$X_0 = f^{-1}(B_0).$ Then it is clear that there is a 2-form $\eta_0$ on $B_0$
such that
$$ \omega \vert B_0 = f^*(\eta _0).$$
Let $\tilde S$ be the closure of the set of all $b \in B$ such that $f^{-1}(b)$ is a multiple
elliptic or a singular rational (possibly reducible) curve. Let $\tilde S_i$ be its irreducible
components.   
Since everything is algebraic, $\eta _0$ has a meromorphic extension to $B.$ 
Since $ B \setminus B_0 = \tilde S \cup E$ with a finite set $E,$ this means that
$\eta \in \tilde \Omega^2_B(\sum k_i\tilde S_i).$ To be more precise, we have first that 
$\eta \vert {\rm reg}B \in H^0({\rm reg}B, \Omega^2_B(\sum k_i \tilde S_i)),$ and then $\eta$ extends
to all of $B,$ since $\tilde \Omega^2_B(\sum k_i \tilde S_i)$ is reflexive. However $\eta$ cannot have
poles along components $T$ of $\tilde S$ which are not in $S$ because the fiber over general 
points in $T$ are generically smooth; therefore the poles cannot cancel.  

\bn {\it Remark.} Actually we have also $\omega \in \Omega^2_B(kS).$ In fact,
let $b \in B$ be a singular point. Then locally around $b = 0,$ we have $B = {\bf C^2}/G$
with a finite group $G \subset Gl(2).$ Let $p: {\bf C^2}Ê\la B$ be the projection. Then
$p^*(\eta \vert B \setminus \{ 0\})$ extends around $p^{-1}(0)$ to a meromorphic form $\hat \eta \in
\Omega^2(k \hat S)$ with $hat S = p^{-1}(S).$ Since $p^*(\eta)$ is $G-$invariant, $\hat \eta$
is $G-$invariant, our claim follows. 

\bn {\bf 5.2 Lemma} {\it Assume that the j-invariant of $f$ is constant. Then there is a
finite set $E \subset B$ such that $f: X_0 \la B_0 = B \setminus E$ is almost smooth.
There is a finite cover $\tilde B \la B,$ such that the pull-back
$\tilde X \la \tilde B$ is an elliptic bundle over $\tilde B_0.$ }

\bn {\bf Proof.} Let $C \subset B$ be an ample divisor and consider $X_C = f^{-1}(C).$ 
Then the j-invariant of $X_C$ is constant so that $f \vert X_C$ is almost smooth ([BPV84]
e.g.).This proves the existence of $E.$ The remaining part is standard. 

\bn {\bf (5.2.a)} If in (5.2) $\omega$ is $f-$vertical, then we conclude by (5.1) that
$\omega = f^*(\eta)$ with a {\it holomorphic} 2-form $\eta$ on $B.$ In the general case
we can say the following; possibly after some cover. Choose a point $b \in B$ outside the finite
set described in (5.2) and let $U$ be an open 
neighborhood such $X_U = f^{-1}(U) \simeq U \times E$ with $E$ the typical fiber of $f.$ 
Let $p_i$ denote the obvious projections. Then
$$ \omega \vert X_U = p_1^*(\omega_1) + p_1^*(\tau_1) \wedge p_2^*(\tau_2)$$
with a 2-form $\omega_1$ on $U$ and 1-forms $\tau_i$ on $U$ and $E,$ respectively. 
This is also described by the exact sequence
$$ 0 \la f^*(\Omega^2_B) \la \Omega^2_X \la \Omega^1_{X \vert B}Ê\otimes f^*(\Omega^1_B) \la 0$$
on $B_0.$ Thus $\omega$ is induced in the sense of (1.2) and Remark (2) after (1.2). 

\bn Concerning non-vertical 2-forms we first show

\bn {\bf 5.3 Lemma} {\it Assume that the j-invariant $J: B \rightharpoonup \bP_1$ of 
$f$ is non-constant. Let $\hat B \la \bP_1$ be a holomorphic model of $J.$ Then there exists 
a surface $g: S \la \bP_1$ and a meromorphic dominant map 
$X \rightharpoonup S \times_{\bP_1}Ê\hat B.$}

\bn {\bf Proof.} To keep notations easy, we assume from the beginning that $J$ is holomorphic (note
that a priori $J$ is almost holomorphic, i.e. proper and holomorphic over an open set of $\bP_1)$. 
Let $g = J \circ f.$ Let $G$ be the general fiber of $g,$ say $G = g^{-1}(c).$ Then we have
an elliptic fibration $\tau_c: G \la J^{-1}(c)$ whose $j-$invariant is constant by construction.
Since $K_X$ is nef, $K_G = K_X \vert G$ is nef, hence $\tau_c$ is relatively minimal and therefore
the only singular fibers if $\tau$ are multiples of smooth elliptic curves. 
\sn We now construct the surface $S$ birationally. Vaguely, we take the collection of the 
general fiber $E_c$ of each $\tau_c$ for general $c$ and then take closure. To be more precise,
we consider the meromorphic map $h: \bP_1 \la \C,$ associating to the general point $c$  the
elliptic curve $E_c$ in the cycle space $\C.$ Let $A$ be the closure of the image of $h$ and
$p: T \la A$ the associated family so that the general fiber of $p$ is of type $E_c.$ After
normalising $A,$ the map $h: \bP_1 \la A$ is holomorphic. 
Now let $S' = T \times_A \bP_1$ and let $S$ be a desingularisation of $S'.$ 
We clearly have a meromorphic dominant map $X \rightharpoonup S \times_{\bP_1}B.$  

\bn {\bf 5.4 Corollary}  {\it In the situation of (5.3) let $\tilde J: \hat B \la C$ be the
Stein factorisation of $J: \hat B \la \bP_1.$ Then $S \la \bP_1$ induces a base change plus 
desingularisation $\tilde S \la C.$ Then there is a birational rational map
$\rho: X \rightharpoonup \tilde S \times _C \hat B.$ We can choose $\tilde S$ such that the
fiber product $\tilde S \times _C \hat B$ is smooth.}

\bn {\bf Proof.} Everything is clear except for the smoothness statement. Here we go back to
the proof of (5.3) and remark that applying a general automorphism of $\bP_1,$ we may
assume that the singular fibers of $S \la \bP_1$ and $J$ are over different points of $\bP_1,$
hence $S \times_{\bP_1} B$ is smooth. Therefore also $\tilde S \times_C B$ is smooth. 

\bn {\bf (5.5)} By (5.4) every elliptic fibration with $\kappa = 2$ and non-constant j-invariant
is birational to a smooth fiber product. Therefore we are reduced to the case of a fiber 
product $X = S_1Ê\times_C S_2$ with maps $g_i: S_i \la C.$ Let $f_i: X \la S_i$  
be the projections. Then every 2-form $\omega$ on $X$ is of the form
$$ \omega = f_1^*(\omega_1) + f_2^*(\omega_2) + f_1^*(\eta_1) \wedge f_2^*(\eta_2)$$
with 2-forms $\omega_i$ on $S_i$ and 1-forms $\eta_i$ on $S_i.$
In fact, we compute $R^2f_{1*}(\cO_X)$ via the Leray spectral sequence:
$$ H^2(\cO_X) = H^2(\cO_{S_2}) \oplus H^1(S_2,R^1f_{1*}(\cO_X)). \eqno (1)$$
Note here that $R^2f_{1*}(\cO_X) = 0$ since all fibers of $f_1$ are 1-dimensional. 
It remains to compute the second term in (1). We apply the Leray spectral sequence for
$g_2$ and use $$R^1f_{1*}(\cO_X) \simeq g_1^*R^1g_{2*}(\cO_{S_2})$$ 
(compare (1.5)). Therefore
$$ H^1(R^1f_{1*}(\cO_X)) = H^1(R^1g_{2*}(\cO_{S_2})) \oplus H^0(R^1g_{1*}(\cO_{S_1}) \otimes
R^1g_{2*}(\cO_{S_2})).$$
Since the first term is clearly $H^2(\cO_{S_2})$ and the second term corresponds to wedge 
products of 1-forms, our claim follows. 

\bn We discuss now the situation when $\omega$ has only finitely many zeroes, i.e.
$Z = \{\omega = 0\}$ is finite. If we  start with a projective smooth threefold $X$ with
$\kappa (X) = 2,$ then by (3.3) $K_X$ is nef, so we again assume now $X$ to be a minimal
model, possibly singular. 

\bn {\bf 5.6 Proposition}  {\it If $\omega$ is vertical, then $f$ has only finitely many 
singular fibers which are not multiple and these singular fibers have only finitely many singularities, all contained
in $Z.$ }

\bn {\bf Proof.} Our claim is local in $B.$ Since $B$ has only quotient singularities, it
is thus clear that we many assume $B$ smooth. Let
$$ T = \{ x \in X \vert {\rm rk}\Omega^1_{X \vert B,x} \geq 2Ê\}.$$
Let $S$ be the union of the multiple fibers. 
Then we must show that $T \setminus S \subset N.$ Take $b \in B$ such that there is some $x \in T \setminus S$ with
$f(x) = b.$ Since $\omega$ is vertical, we have $\omega = f^*(\eta)$ with $\eta \in H^0(\omega_B(\sum k_i S_i))$
where the $S_i$ are the irreducible components of the closure of $S.$  First suppose that $x$ is not contained in the
closure of $S.$ Then  locally we can write
$$ \eta = \sum \alpha_i \wedge \beta_i$$ 
with holomorphic 1-forms $\alpha_i$ and $\beta_i.$ 
Since the kernel of $$\Omega_X^1 \vert \{x\} \la \Omega^1_{X \vert B} \vert \{x\} $$ is 1-dimensional, we conclude
$$F^*(\alpha_i)(x) \wedge f^*(\beta_i)(x) = 0,$$ what had to be proved. 
If $x \in {\overline S},$ then we perform a (local) base change in $B$ to kill the multiplicities and apply
the previous argument on $\tilde XÊ\la \tilde B.$

\bn \bn \bn {\medium 6. Torus fibrations}
\bn {\bf (6.1)} Let $r : X \to S$ be a surjective holomorphic map with connected fibers
between
compact complex connected manifolds. We assume throughout this section that $X$ is KŠhler with 
generic
fiber $F = F_s = X_s := r^{-1}(s)$ a complex torus.

\noi For the relevant notions on relative cycle (or Douady) spaces we refer
to [Ca85] and
[Fu83] where the constructions used here have been introduced with more
details.
\bn {\bf (6.2)} Let $r_\ast : \hbox{Aut}_S(X) \to S$ be the relative
automorphism space of
$S$ over $S$; it consists of the closure in ${\cal C}\Big(X \dis\mathop
{\times}_{S}
X/S\Big)$ of the graphs of automorphisms $\alpha_S : F_s \to F_s$
of some
smooth fiber of $r$, where these graphs are considered in ${\cal C}(F_s \times F_s) \subset
{\cal C}\Big(X
\dis\mathop {\times}_{S} X/S\Big)$ of automorphisms $\alpha_S : F_S \to F_S$
of some
smooth fiber of $r$.
\smallskip
\noi By the KŠhler assumption on $X$, Aut$_S(X)$ has compact irreducible
components
([Li78]). There is one distinguished component $r'_0 : X'_0 \to S$
of Aut$_S(X)$,
which lies surjectively on $S$ and has a canonical meromorphic zero-section
$\zeta'_0 : S
\to X'_0$, namely the component containing id$_{F_S}$ for $s$ generic in $S$. We
denote by $r_0 :
X_0 \to S$ (and $\zeta_0 : S \to X_0$) any smooth model of $X'_0$,
isomorphic to $X'_0$
over some non-empty Zariski open subset $S^\ast$ of $S$.

\bigskip \noindent {\bf 6.3 Definition} We call $r_0 : X_0 \to S$ (equipped with
$\zeta_0$) the
{\it untwisted model of} $X$.

\bn  {\bf 6.4 Proposition} {\it Let $S'$ be a Zariski dense open subset
of $S$ over which
$r$ (and hence $r_0$) is smooth. Let $r' : X' \to S'$ be the restriction of
$r$ over $S'$ to $X'
:= r^{-1}(S')$ ; use similar notations for $r_0$.\hfill\break
For $p \ge 0$, there are canonical isomorphisms $$\psi'_p: r'_\ast
(\Omega^p_{X'/S'}) \to
(r'_0)_\ast (\Omega^p_{X'_0/S'}),$$ and by Hodge symmetry: $$\varphi'_p : R^p
r'_\ast ({\cal
O}_{X'/S'}) \to (R^p r'_0)_\ast ({\cal O}_{X'_0/S'})$$.}

\noi {\bf Proof.} The assertion is of local nature in $S'$. It is
sufficient to check it over
any sufficiently small open neighborhood $S''$ of any $s \in S'$. But then
one has
isomorphisms $\psi : r^{-1}(S'') \to r^{-1}_0(S'')$, unique up to
translations in the fibers,
which define the asserted isomorphisms.

\bn {\bf 6.5 Theorem} \  {\it Assume that $r$ and $r_0$ are locally
projective. Then 
$h^p({\cal O}_X) = h^p({\cal O}_{X_0})$ for all $p \ge 0$ and in particular $\chi({\cal O}_X) =
\chi({\cal
O}_{X_0})$.}

\noi {\bf Proof.} We can assume that $S-S'$ is a divisor with only normal
crossings. Then the
sheaves $R^p r_\ast {\cal O}_{X/S}$ are locally free and coincide with
the canonical
extensions of their restrictions to $S'$ ([Ko86], [Mw87], [Na86];
observe that the
assertion is of local nature over $S$).

\noi The same properties hold for $r_0$, so that $\varphi'_p$ extend to
isomorphisms
$\varphi_p$. These isomorphisms define naturally an isomorphism of the Leray spectral sequences 
for ${\cal O}_X$ and ${\cal O}_{X_0}$ with respect to the maps $r$ and $r_0.$ Hence 
$$H^\ast({\cal O}_X) \simeq 
H^\ast({\cal
O}_{X_0}),$$ establishing our claim.
\medskip
\noi {\bf Remark.} The assertions of (6.5) hold true certainly in the KŠhler
case as well;
unfortunately the necessary tools seem to have been only written up in the
projective case.

\bn {\bf 6.6 Corollary} {\it Assume that the generic fiber of $r$ is an
abelian variety.
Then the conclusions of 6.5 hold.}

\bn {\bf Proof.} Since $X$ and $X_0$ are K\"ahler and $F$ is abelian, $r$ and $r_0$ are locally
Moishezon over
$S$. Hence suitable bimeromorphic models are locally projective over $S$
and 6.5 applies.

\bn {\bf 6.7 Proposition} {\it Let $r : X \to S$ be as in 6.1. Assume
moreover that the
smooth fibers of $r$ are all isomorphic to a {\it fixed} torus $T$. Then
:\hfill\break
{\rm (i)} There exists a finite Galois cover $\beta :
\widetilde S \to S =
\widetilde S/G$ with Galois group $G$ such that the induced map $\tilde r: \widetilde X := X \dis\mathop
{\times}_{S}
\widetilde S \to \widetilde S$ is bimeromorphically a principal fiber bundle with fiber $T$.\hfill\break
{\rm (ii)} Let $\beta$ be as in (a) and $\tilde r_0 := \widetilde X_0 \to
\widetilde S$ be the
untwisted model of $\tilde r$ (also obtained from $r_0$ by the base change $\beta$). Then $\tilde r_0$ is
bimeromorphic to the second projection $p_2 : T_0 \times \widetilde S \to \widetilde S$ where $T_0 = \hbox{\rm Aut}^0(T)$.
\smallskip
\noi {\rm (iii)} $h^0(\widetilde X,\Omega^p) = \dis\sum_{q+t=p} h^0(\widetilde
S,\Omega^q) \times h^0(T,\Omega^t) = h^0(\widetilde X_0,\Omega^p)$ $(p \ge
0)$.\hfill\break
{\rm (iv)} $h^0(X,\Omega^p) = h^0(X_0,\Omega^p)$ $(p \ge 0)$.}

\bn {\bf Proof.} (i) Let $$r_\ast : \overline J := \hbox{Isom}_S(T \times
S,X) \to S$$ be the
space of relative isomorphisms over $S$ of $T \times S$ and $X$. Over a
smooth fiber $F$
of $r$ it consists of all isomorphisms of $T$ and $F$, each identified
with its graph and
considered as a point in $${\cal C}(F \times T) \subset {\cal C}\Big((T
\times S) \dis\mathop
{\times}_{S} X = T \times X\Big).$$ Again see [Ca85] or [Fu83] for
further details. Let
$r_\ast : J \to S$ be some irreducible component of $\overline J$ lying
surjectively over $S$
(such a component exists by the KŠhler assumption on $X$ and the fact that
$F_s = F$ is
isomorphic to $T$ for $s$ generic in $S$).

\noi Let $\rho : J \to \overline S$ be the Stein factorisation of $r_\ast :
J \to S$. We have
a natural meromorphic composition map $\tau : T_0 \times J \to J$ over $\overline S$, which
generically over $\overline S$ maps $(j_s : T \to F_s) \in
\rho^{-1}(\overline s)$ to $\tau(t,j_s) := j_s \circ t.$ Here $t
\in T_0 =
\hbox{Aut}^0(T)$.

\noi The conclusion now follows from the facts that (as just shown)  $\rho
: J \to
\overline S$ is bimeromorphically a principal fiber bundle with fiber $T_0$, and that the
natural evaluation
map over $\overline S$ at $t'$ $\varepsilon : J \to \overline X := \Big(X
\dis\mathop
{\times}_{S} \overline S\Big)$ is an isomorphism (bimeromorphically), for
any fixed $t' \in
T$, where $\varepsilon$ sends $j_s$ to $j_s(t')$.

\noi We conclude by taking a Galois base change dominating $\overline S
\to S$.

\noi (ii) This immediately follows from the fact that $\tilde r_0 :
\widetilde X_0 := X
\dis\mathop {\times}_{S} \widetilde S \to \widetilde S$ is the untwisted
model of $\tilde
r : X \dis\mathop {\times}_{S} \widetilde S \to \widetilde S$, and is thus
bimeromorphically a principal
bundle of fiber $T_0$ with a (meromorphic) section.

\noi (iii) This follows from the results of A. Blanchard ([Bl58]) about the
degeneration of the
Leray spectral sequence.

\noi (iv) Observe that we do not assume here $r$ to be locally projective, otherwise (iv)
follows from (6.5). The Galois group $G$ of the cover $\beta : \widetilde S \to S'$
acts in a natural
way on $H^0(T,\Omega^p)$ $(p \ge 0)$ as follows : let $\overline s \in
\overline S$ and $g
\in G$. Choose $j_{\overline s} : T \to F_S$ and $j_{\gamma\overline s} : T \to
F_S$ arbitrarily
and define $\tilde g := (j^{-1}_{\gamma\overline s} \circ j_{\overline s}) \in
\hbox{Aut}(T)$. $\tilde g$ is well defined up to a translation in $T$. Its action $g^\ast$ on
$H^0(F,\Omega^p)$ is thus well-defined, independent on the choices $(\widetilde
S,j_{\tilde s},j_{\gamma\tilde s})$ made, as one checks easily. Observe that this
action is the same for $\widetilde X$ and for $\widetilde X_0$. We then
just deduce the
claim (iv) from the fact that $H^0(X,\Omega^p) = H^0(\widetilde
X,\Omega^p)^G$, where
the action is defined via the decomposition $H^0(\widetilde X,\Omega^p) =
\dis\bigoplus_{q+t=r} H^0(\widetilde S,\Omega^q) \otimes H^0(T,\Omega^t)$, and
similarly for $X_0$.
\bn {\bf  Remark.} The equality $q(\widetilde X) = q(\widetilde S) +
q(T)$ was
obtained in [Fu83] using Deligne's mixed Hodge structures in the K\"ahler
situation. We
gave here a more elementary proof, avoiding the K\"ahler version of the
theory of Deligne
(which we were not able to find explicitely written in the literature).

\bn {\bf 6.8 Corollary} {\it Let $r : X \to S$ be a surjective connected
holomorphic map
from a KŠhler threefold $X$ to a curve $S$.\hfill\break
Assume that its generic fiber $F$ is either a torus or a K3-surface and
that $a(F) \le
1$.\hfill\break
Then any two smooth fibers of $r$ are isomorphic.}

\bn {\bf Proof.} We follow here an argument of [Fu83].\hfill\break
Assume first that $F$ is a K3-surface. Since $X$ is not projective, we have
$h^{0,2}(X) > 0$. If
$0 \not= \omega \in H^0(\Omega^2_X)$, then $\omega|_{\dis F} \not= 0$
since $q(F) =
0$, and $d\omega = 0$, since $X$ is K\"ahler. Hence the periods of $\omega$
on $F = F_s$
are constant in a local marking, and the conclusion follows from the
Torelli theorem for
K3-surfaces.\hfill\break
If $F$ is a torus, we replace $X$ by $X_0$ and consider the ${\Bbb
Z}_2$-quotient $p' :
Y'_0 := (X_0/\pm 1) \to S$ of $r_0 : X_0 \to S$ by the Kummer involution in
the fibers of
$r_0$, relative to the canonical section (an easy cycle space argument
shows that this
involution extends meromorphically over all of $S$). Then replace $p'$
by a smooth
model $p : Y_0 \to S$ with generic fiber $Y_{0,s}$, the Kummer surface of
$X_{0,s}$.
Because $a(F) \le 1$, we still have $h^{0,2}(Y_0) > 0$ ; let $0 \not= \omega \in
H^{0,2}(Y_0)$ and let $\overline \omega$ be its lift to $X_0$. Then argue
as above, using
the Torelli theorem for 2-forms on 2-dimensional complex tori.
\medskip
\noi We now consider the Kodaira dimension of $X_0$ :

\bn {\bf 6.9 Proposition} {\it Let $r : X \to S$ be as in 6.1 with $\dim X = \dim S +1,$ i.e. $r$ is
an elliptic fibration. Then $\kappa(X_0) \le \kappa(X)$.}

\bn {\bf Proof.} By [Ue87],[Na88] we have a canonical bundle formula :
$$K_X = r^\ast(K_S \otimes r_\ast \omega_{X/S}) \otimes {\cal
O}_X\big(\Sigma(m_i-1)
D_i\big) \otimes {\cal O}_X(\Delta),$$
where :\hfill\break
$\bullet$ $\Delta$ is an effective divisor supported on the 2-dimensional
fibers of
$r$.\hfill\break
$\bullet$ $r_\ast(\omega_{X/S})$ is locally free (assuming $S-S'$ to have only normal
crossings).\hfill\break
$\bullet$ The divisors $D_i$ are those having generically over their image
$\Delta_i$ in
$S$ multiple fibers of multiplicity $m_i \ge 2$.\hfill\break
Thus : $\kappa(X) = \kappa(S,D)$ where $D$ is the ${\Bbb Q}$-divisor
defined by : $D = K_S
\otimes r_\ast(\omega_{X/S}) \otimes {\cal O}_S\bigg(\Sigma\Big(\dis{m_i-1 \over
m_i}\Big) \Delta_i\bigg)$.\hfill\break
We perform the same arguments on $X_0.$ Since $r_0$ has a meromorphic section, it
has no multiple fibers, and thus we get
$$ \kappa X_0) = \kappa(S,K_S \otimes (r_0)_*(\omega_{X_0 \vert S}).$$
Since $r_*(\omega_{X \vert S}) \simeq (r_0)_*(\omega_{X_0 \vert S})$ by (6.4), our
claim follows. 
\bn {\bf  Remarks.} (i) One may have strict inequality: if $X$
is an Enriques
surface, then $X_0$ is a rational surface (a classical fact). Indeed: we have $q(X) = q(X_0) =
0$ and $\kappa(X_0) = -\infty$ since deg$(r_\ast
\omega_{X/S}) = \chi({\cal O}_S) = 1$.\hfill\break
(ii) The inequality 6.8 is probably true in general, i.e. for higherdimensional fibers, but we are able to
show this only in
special cases.

\bn {\bf 6.10  Lemma} {\it Let $X$ be a normal compact threefold and $h :
X \to {\Bbb
P}_1$ and $g : X \to S$ be connected holomorphic surjective maps. Assume
that $S$ is a
normal surface with $h^0(K_S) \not= 0$, and that {\it every} fiber $F$ of $h$ is a
finite quotient of either a torus or a K3-surface. Assume moreover that $h$
has a
section and that $g(F) = S.$\hfill\break
Then $h \times g : X \to {\Bbb P}_1 \times S$ is bimeromorphic.}

\bn {\bf Proof.} It will be sufficient to show that if $G$ is a generic
fiber of $g$, then
$h|_{\dis G} : G \to {\Bbb P}_1$ has degree $d$ equal to one.\hfill\break
Let $0 \not= \omega \in H^0(K_S)$ and $\eta = g^\ast (\omega)$ is at least a 2-form
on the regular part of $X$ but it defines also restrictions $\eta |Ê{\dis F} \in
H^0(K_F.)$ Moreover $\eta $ defines a section 
$$ s \in H^0(K_{X \vert {\Bbb P}}).$$ 
\hfill\break
For the general reduced fiber $F$ of $h$, $\eta |_{\dis F}$ is everywhere non-zero ( as 2-form on the
smooth part
of $F$ or as section of $K_F$). On the other hand, we have
$$(K_{X/{\Bbb P}_1})\cdot{\dis G} = 2\gamma -2 + 2d$$ where $\gamma$ is the genus
of $G$,
because :$K_{X/{\Bbb P}_1}  \simeq  K_X \otimes h^\ast {\cal
O}_{{\Bbb P}_1}(2)$. If $d \ge 2$, then $s|_{\dis G}$ has zeroes, so that
$s|_{\dis F}$ has zeroes, hence vanishes identically for general $F.$ This is a
contradiction.

\vfill \eject \noindent  {\medium 7. Algebraic reductions of K\"ahler threefolds}
\bn  We shall collect in this section some facts about algebraic reductions
of compact
K\"ahler manifolds, in particular in dimension 3. In the threefold case, most are due to
[Fu83]; we will
provide here short and elementary proofs.

\bn{\bf 7.1 Theorem {\rm ([Ue75])}.} {\it Let $r : X \to S$ be the
algebraic reduction of a
compact complex manifold, and $F = X_s$ a general fiber of $r$. Then 
$\kappa(F) \le 0$.}

\bn {\bf 7.2 Theorem} {\rm ([Ca85], [Fu83] when $Y$ is
projective)} {\it Let $f : X \to
Y$ be a surjective connected map between compact complex manifolds with $X$
K\"ahler.
Assume that the generic fiber $F$ of $f$ is projective. Then :\hfill\break
{\rm (i)} If $q(F) = 0$, then $a(X) = a(Y) + {\dim} F$.\hfill\break
{\rm (ii)} If $a(X) = a(Y)$, then $F$ is almost-homogeneous.\hfill\break
{\rm (iii)} If $X$ contains a compact analytic set $Z$ with $\dim Z = \dim Y$
and $f(Z) = Y$,
then $a(X) = a(Y) + {\dim} F$.}

\bn {\bf 7.3 Corollary} {\it If $r : X \to S$ is an algebraic reduction
of $X$. If its
generic fiber $F$ is projective, then we have:\hfill\break
{\rm (i)} $q(F) > 0$ ; $F$ is almost-homogeneous and there does not exist
any $Z \subset X$ dominating $S.$ 
\hfill\break
{\rm (ii)} In particular, if $\dim F = 2$, then $F$ is either an abelian
surface or a ruled
surface over an elliptic curve $E$, isomorphic to ${\Bbb P}({\cal O}_E
\otimes {\cal L})$,
where ${\cal L} \in \hbox{\rm Pic}^0(E)$ is not torsion.}

\bn {\bf Proof.} (7.3) follows directly from (7.2) and the classification of
surfaces.

\bn {\bf 7.4 Corollary.} {\it Let $r : X \to S$ be the algebraic
reduction of a compact
K\"ahler threefold with $a(X) = 1$. Then the generic fiber $F$ of $r$ is
either \hfill\break
{\rm 1.} bimeromorphic to a torus or a K3-surface of algebraic
dimension
zero,\hfill\break
{\rm 2.} a ruled surface ${\Bbb P}({\cal O}_E \otimes {\cal L})$ as in
7.3.(ii).}

\bn {\bf Proof.} This follows from the classification of surfaces with
$\kappa \le 0$
again : Enriques and rational surfaces are excluded because projective with
$q = 0$ ;
bielliptic (hyperelliptic) surfaces do not occur because they are projective but not almost-homogeneous. The
almost-homogeneous birationally ruled surfaces with $q > 0$ are exactly those of
7.3.(ii).\hfill\break
Finally if $a(F) \le 1$, it has to be bimeromorphic to either a torus or a
K3-surface. In this
last case (K3), $a(F) = 1$ is excluded, because taking a relative algebraic
reduction ([Ca81]) we obtain a meromorphic factorisation $X \to S \to C$ with
$S \to C$ being a ${\Bbb P}_1-$fibration, so that $a(S) = 2$ contradicting $a(X) = 1.$ 

\bn The special case 7.4.(i) has been analyzed more closely in [Fu83] :

\bn {\bf 7.5 Lemma} {\rm ([Fu83])} {\it Let $r : X \to S$ be the algebraic reduction of a
compact K\"ahler threefold with $a(X) = 1.$ 
and assume that
$F$ is bimeromorphic to either a torus or a K3-surface. Then $F$ is either
a torus or a
K3-surface (after taking a suitable bimeromorphic model of $X$).}

\bn {\bf Proof.} In the torus case, we just replace $r : X \to S$ by its
Albanese reduction ([Ca85]).\hfill\break
In the K3-surface case, let $0 \not= \omega \in H^{0,2}(X)$. By the
arguments used to
prove 6.8, we see that any smooth $F_s$ is bimeromorphic to a fixed K3-surface
$F_0$.\hfill\break
Consider now any irreducible component $\mu : M \to S$ of the space of $S-$morphisms
Mor$_S(F_0
\times S,X) \subset {\cal C}(F_0 \times X)$ which lies over $S$ and consists generically over
$S$ of graphs
of bimeromorphic maps from $F_0$ to $F_s = X_s$. $\mu : M \to S$ be the
canonical map.\hfill\break
Then $\mu : M \to S$ is onto and generically finite, and there is a natural
evaluation map
$\varepsilon : F_0 \times M \to \widetilde X := X \dis \mathop
{\times}_{S} M$,
which is bimeromorphic. By a further finite base change, we can assume that
$\mu : M
\to S$ is Galois cover with Galois group $G$.\hfill\break
We can now define, using the construction of the proof of 6.7.(iv), an
action of $G$ on
$F_0 \times M$ such that $F_0 \times M/G$ is bimeromorphic to $X$ over
$S = M/G$.

\bn {\bf 7.6 Corollary} {\rm ([Fu83])} {\it Let $X$ be a non-projective
compact K\"ahler
threefold. Then $X$ is bimeromorphic to one of the following, according to
its algebraic
dimension $a(X)$ :\hfill\break
{\rm (i)} $\underline{a(X) = 0}$ {\rm : (i.a)} simple non-Kummer.\hfill\break
$\hphantom {''''''.................}$ {\rm (i.b)} simple Kummer.\hfill\break
$\hphantom {''''''.................}$ {\rm (i.c)} a holomorphic fibration $f : X \to
\Sigma$ to a normal surface with
$a(\Sigma) = 0$ and generic fiber ${\Bbb P}_1$.\hfill\break
$\hphantom {''''''.................}$ This fibration is unique.\hfill\break
Let $F$ be a generic fiber of $r : X \to S$, the algebraic reduction of
$X$.\hfill\break
{\rm (ii)} $\underline{a(X) = 1}$ {\rm : (ii.a)} $a(F) = 0$ and $F$ is a
fixed torus.\hfill\break
$\hphantom {''''''''i...............}$ {\rm (ii.b)} $a(F) = 0$ and $F$ is a
fixed
K3-surface.\hfill\break
$\hphantom {''''''''i...............}$ Then $X = F\times \widetilde S/G$
for some Galois cover
$\widetilde S$ of $S = \widetilde S/G$.\hfill\break
$\hphantom {''''''a................}$ {\rm (ii.c)} $a(F) = 1$ and $F$ is a
fixed torus.\hfill\break
$\hphantom {''''''''i...............}$ {\rm (ii.d)} $a(F) = 2$ and $F$ is a
(possibly varying)
abelian surface.\hfill\break
$\hphantom {''''''''i...............}$ {\rm (ii.e)} $a(F) = 2$ and $F =
{\Bbb P}({\cal O}_E \oplus {\cal
L})$, as in 7.3.(ii).\hfill\break
{\rm (iii)} $\underline{a(X) = 2}$ {\rm :} $F$ is an elliptic curve
(possibly varying).}

\bn {\bf Proof.} Everything has been proved $\big($using 6.8 for $a(X) =
1$, $a(F) = 0,1$,
and additionally the proof of 7.4 if $F$ is K3 with $a(F)=0\big)$, except
that if $a(X) = 0$
and $X$ is not simple, then there exists $f : X \to \Sigma$ as in (i.c).
Because $X$ is not
simple, it is covered by a family of curves. Only one curve goes through
the generic point
of $X$, otherwise $X$ were covered by divisors, contrary to $a(X) = 0$. The
curves on $X$
thus define a (unique) meromorphic map $f : X \to \Sigma$ to a surface
$\Sigma$ with $a(\Sigma) = 0$, since $0 = a(X) \ge a(\Sigma)$. Let $g$ be
the genus of a
generic fiber of $f$ ; because $a(\Sigma) = 0$, any two smooth fibers of
$f$ are isomorphic
to a fixed curve $C$. If $g \ge 2$, then there exists a finite cover $\widetilde \Sigma \to
\Sigma$ such that $\widetilde X := \Big(X \dis \mathop
{\times}_{\Sigma} \widetilde \Sigma\Big) = \widetilde \Sigma \times C$,
since Aut$(C)$ is
finite. This contradicts $a(X) = a(\widetilde X) = 0$. If $g=1$, first notice that $f$
cannot have moduli, since otherwise $\kappa (X) > 0$ by $C_{2,1}$ contradicting $a(X) = 0.$
Hence $q(\widetilde X) =
q(\widetilde \Sigma) + 1$ by 6.7, and thus $q(\widetilde \Sigma) = 0$ or 2 since
$a(\widetilde
\Sigma) = a(\Sigma) = 0$, i.e. bimeromorphically a torus or K3. If $q(\widetilde \Sigma) = 0$, we get
$q(\widetilde X) = 1$
contradicting $a(\widetilde X) = 0\;;$ if $q(\widetilde \Sigma) = 2$, then
$\widetilde X$ is
bimeromorphic to a torus since $a(\widetilde X) = 0$, and $X$ is Kummer
contradicting
our assumption. Thus $g = 0$ and $C = {\Bbb P}_1$, as claimed. The fact that we can
choose $f$ holomorphic is easily obtained by contracting some curves in $\Sigma.$ 

\vskip 1cm

\noi {\medium 8. K\"ahler threefolds with $\kappa = 0,1$.}
\bn {\bf (8.0) Convention}  $X$ will always denote a compact K\"ahler connected manifold with
$\dim X= 3$ and we assume that $X$ is not both simple and not Kummer. Let
$r : X \to S$
be an algebraic reduction of $X$ and $r_0 : X_0 \to S$ be its untwisted model. 
Both maps may assumed to be holomorphic.

\bn  {\bf 8.1 Theorem} {\it Let $X$ be a compact K\"ahler
threefold which is not both simple and non-Kummer. Then :\hfill\break
{\rm (i)} If $\kappa(X) = -\infty$, $X$ is uniruled.\hfill\break
{\rm (ii)} If $\kappa(X) = 0$ and $h^{0,2}(X) \not= 0$, then $X$ is
bimeromorphic to some threefold
$X'$ (possibly with quotient singularities) which has a finite cover
$\widetilde X'$ Žtale in
codimension one such that $\widetilde X'$ either a torus or a product of
an elliptic curve
by a K3-surface.}

\noi {\bf Proof.} Because the results are known when $X$ is projective, due
to [Mi88],
[Mo87], [Ka85] and [Vi80], we consider only the non-algebraic case and proceed
according to
$a(X) = 0, 2, 1$ successively in the increasing order of difficulty.

\bn  {\bf (8.2)} $\underline{a(X) = 0}$. Then by (7.6) $X$ is either Kummer with
$\kappa(X) = 0$ so
that 8.1.(ii) holds, or $X$ is a
generic ${\Bbb P}_1$-bundle over a surface $\Sigma$, and 8.1.(i) holds.

\bn {\bf (8.3)} $\underline{a(X) = 2}$. Notice then that $\big($by
7.2.(iii)$\big)$ $X_0$ is
projective since $r_0$ has a section. Moreover, $\kappa(X_0) \le \kappa(X)$
by (6.9) and
$h^{0,2}(X_0) = h^{0,2}(X) > 0$ (by 6.5 and the assumption that $X$ is not projective).

\noi We thus distinguish the following possible cases :

\noi $\kappa(X) = \kappa(X_0) = -\infty$ (8.4), \hfill \break 
$\kappa(X) = 0 >
\kappa(X_0) = -\infty$
(8.5) \hfill \break
and $\kappa(X) = \kappa(X_0) = 0$ (8.6).

\bn {\bf (8.4)} $\underline{a(X) = 2 \, ; \kappa(X) = \kappa(X_0) =
-\infty}$. Then $X_0$ is uniruled, $X$ being projective.

\noi Let $\sigma_0 : X_0 \to \Sigma_0$ be its rational quotient ([Ca92],
[KoMiMo92]). Because $h^{0,2}(X_0) = h^{0,2}(X) > 0$, $X$ being non-algebraic, $\Sigma_0$ 
is a surface with
$h^{0,2}(\Sigma_0) > 0$,
and the generic fiber of $\sigma_0$ is ${\Bbb P}_1$. We may assume $\Sigma_0$ smooth and
$\sigma_0$ holomorphic.

\noi Let $F_0$ be the generic fiber of $r_0$, and $E := \sigma_0(F_0)$. Then $E$ 
is a member of a
1-dimensional family of elliptic curves on $\Sigma_0$. Notice that $E$ is
elliptic and not
${\Bbb P}_1$ because $\kappa(\Sigma_0) \ge 0\big)$.
Since
$\kappa(\Sigma_0) \ge 0$,
there exists a connected meromorphic surjective map $\tau_0 : \Sigma_0 \to
C$ to some
curve $C$ contracting all the elliptic curves $E.$ In fact, one can introduce an equivalence
relation such that two general points of $\Sigma_0$ are equivalent if and only if they can
bejoined by a chain of these elliptic curves, see [Ca81]. Then not all points are equivalent,
otherwise $\kappa (\Sigma_0) = - \infty$. Hence the meromorphic map $\tau_0$ is just given
by the quotient of the equivalence relation [Ca81].
 
\noi We conclude the existence of a unique map $\rho : S \to C$ such that $\tau_0 \circ \sigma_0
= \rho \circ r_0$.

\noi We thus have constructed a commutative diagram of surjective connected
maps (which may assume to be holomorphic) :
$$\diagram{
&S &\cr
r_0\nearrow & &\searrow\rho \cr
X_0  &\build\hbox to 8mm {\rightarrowfill}_{}^{f_0} &C \cr
\sigma_0\searrow & &\nearrow \tau_0 \cr
&\Sigma_0 & \cr} \leqno \hbox{\bf (8.4.i)}$$
Let $G_0$ be a generic fiber of the composed map $f_0 := \rho \circ r_0$. We have
$q(G_0) =
1$ (because $\tau_0$ is an elliptic fibration and $\sigma_0$ a ${\Bbb
P}_1$-fibration);
also observe that $r_0 : G_0 \to P := r_0(G_0) \subset S$ is the untwisted
model of $r: G
:= r^{-1}(P) \to P$, where $P$ is the generic fiber of $\rho.$

\noi Hence $q(G) = q(G_0) = 1$ by (6.5). Let now $ X \buildrel \sigma \over {\to}
\Sigma \buildrel \tau \over {\to} C$ be the  
Albanese reduction of $f :=
\rho \circ r$. We thus
get a second commutative diagram :
$$\diagram{
&S &\cr
r\nearrow & &\searrow\rho \cr
X  &\build\hbox to 8mm {\rightarrowfill}_{}^{f} &C \cr
\sigma\searrow & &\nearrow \tau \cr
&\Sigma & \cr} \leqno \hbox{\bf (8.4.ii)}$$
Here $\Sigma$ is not algebraic, otherwise $S \dis \mathop {\times}_{C}
\Sigma$
is algebraic and then $X$ would be projective. Thus $h^{0,2}(\Sigma) > 0$, and
$\kappa(\Sigma)
\ge 0$. We then apply $C_{3,2}$ ([Ue87]) to conclude that the generic fiber of
$\Sigma$ has
$\kappa = -\infty$, and is thus ${\Bbb P}_1$. So  $X$ is uniruled, as
claimed.

\medskip

\noi {\bf (8.5)} $\underline{a(X) = 2 \, ; \kappa(X) = 0 \, ; \kappa(X_0) =
-\infty}$. We shall
exclude this case. 
\smallskip \noi Considering $r_0 : X_0 \to S$  and the same arguments as
in 8.4, we get
the diagrams (8.4.i) and (8.4.ii), with the same properties. From $C_{3,2}$ we
get that
$\kappa(\Sigma) = 0$, and that $\sigma : X \to \Sigma$ and $\tau : \Sigma
\to C$ are
elliptic fibrations. Moreover, the generic fiber of $\sigma_0$ is ${\Bbb
P}_1$ and $\tau_0$
is an elliptic fibration since $h^{0,2}(X_0) = h^{0,2}(X) > 0$, so that the
rational quotient
$\Sigma_0$ of $X_0$ is a surface with $h^{0,2} > 0$.

\noi As above in (8.4), the generic fiber $G_0$ of $f_0$ has the following properties: $q(G_0)
= 1$ and
$\kappa(G_0) = -\infty$. The generic fiber of $\rho$ is therefore ${\Bbb P}_1$
since it is the
image of the generic fiber of $\sigma_0$ by $r_0$. Notice that $\kappa (G)  \geq 0$ since
$X$ is not uniruled and that $G$ is an elliptic fibration over ${\Bbb P}_1.$ So $\kappa (G) \leq 1.$
We show that
$\kappa(G) = 0$. Indeed: $G_0$ is bimeromorphic to a ruled surface over an
elliptic curve
$E_0$. This implies that $r_0 : G_0 \to {\Bbb P}_1 = r_0(G_0) \subset S$ is
a generically
(over ${\Bbb P}_1$) locally trivial bundle with fiber $E'_0$ isogeneous to
$E_0$; in
particular, it is a principal fiber bundle over ${\Bbb P}_1$ with fiber
$E'_0$. But $r_0$ is
the untwisted model of $r : G \to {\Bbb P}_1$; the corresponding fiber of
$f$ which is
therefore also a principal fiber bundle with fiber $E'_0$. Hence $\kappa(G) =
-\infty$, because
$\chi({\cal O}_{G_0}) = \chi({\cal O}_G) = 0$. This is a contradiction,
since $\kappa(X) = 0$.

\medskip

\noi {\bf (8.6)} $\underline{a(X) = 2 \, ; \kappa(X) = \kappa(X_0) = 0}$.
Consider again the untwisted model $r_0 : X_0
\to S$. We shall establish the claim successively in the cases
$\kappa(S) \ge 0$;
$\kappa(S) = - \infty$,
$q(S) > 0;$ and finally exclude the case $\kappa(S) = -\infty$, $q(S) = 0$.

\bn (1) Assume first that $\kappa(S) \ge 0$, so that $\kappa(S) = 0$ and $S$
is (up to cover
and bimeromorphy) either a K3-surface or an abelian surface. Moreover $X_0$  is
projective since $r_0$ has a section. Let $r'_0 : X'_0 \to S$ be a minimal model of $X_0.$
Since $X'_0$ has a 2-form $\eta,$ there exists by [Pe94] a 1-form $\omega$ such that
$\eta \wedge \omega \ne 0$ generically. In particular $q(X'_0) > 0.$ Now by [Ka85] there exist a cover, Žtale in
codimension one, say
$$\tilde r'_0 : \widetilde X'_0 = X'_0 \dis \mathop {\times}_{S}
\widetilde S \to
\widetilde S$$ which is a product $\widetilde X'_0 = E \times \widetilde S$
for some
elliptic curve $E$, and $\widetilde S$ either K3 or abelian. But
$\widetilde X'_0$ is the
untwisted model of $\tilde r : \widetilde X := X \dis \mathop {\times}_{S}
\widetilde S \to
\widetilde S$ (see the argument, exposed in the case $\kappa(S) = -\infty$
and $q(S) = 0$
below). This shows the claim in that case.

\bn (2) Let us now assume that $\kappa(S) = -\infty$ and $q(S) > 0$.

\noi Let $\alpha : S \to E$ be the Albanese map for $S$. Because $\kappa(X)
= 0$, we get
from $C_{3,1}$ that $E$ is an elliptic curve. Note that $S$ is (bimeromorphically)
ruled over $E$.
Again from
$C_{3,1}$ we deduce that $F$, the generic fiber of $f = \alpha \circ r : X \to E,$ has
$\kappa(F) = 0$. $F$ has an elliptic fibration over ${\Bbb P}_1.$ Hence $q(F) \not = 2$.
If $q(F) = 1,$ $F$ is bimeromorphically a hyperelliptic surface. Then $q(F_0) = 1$
by (6.4) so that $F_0$ is hyperelliptic or a ruled surface over an elliptic curve. Since
$F_0$ has a section, the first alternative cannot occur. But then $X_0$ is uniruled,
contradicting $\kappa (X_0) = 0.$ So the case $q(F) = 1$ does not occur.

If $q(F) = 0,$ then $F$ is K3 or Enriques. In the latter case $h^q({\cal O}_F) = 0$ for
$q = 1,2,$ so that $R^qf_*({\cal O}_X) = 0$ for $q = 1,2$ (both sheaves a priori being locally free).
Hence $H^2({\cal O}_X) = 0,$ contradition.
If $F$ is K3, we can again apply [Ka85] as above to get our claim.

\bn (3) We now exclude the remaining case $\kappa(S) = -\infty$ and $q(S) = 0$
(i.e. $S$ is
rational). Assume that $S$ is rational; then, arguing as in (1), we have, after changing
$X_0$ birationally, a finite cover, \'etale in codimension 2, say $\widetilde X_0
\buildrel g \over
{\to} X_0$ such that $\widetilde X_0$ is 
a product
$\widetilde E \times \widetilde \Sigma$ where $\widetilde E$ is an elliptic
curve and the
surface $\widetilde \Sigma$ is either K3 or abelian. If we do not insist that $\widetilde X_0$
is really a product, i.e. if we blow up the product, then we may assume that $g$ is holomorphic. 
Let $\tilde r_0 :
\widetilde X_0 \to
\widetilde S \buildrel \sigma \over {\to} S$ with $\sigma \circ \tilde r_0 =
r_0 \circ g$ be the
Stein factorisation of $r_0\chi.$ Let $p_1 : \widetilde X_0 \to
\widetilde E$, $p_2 :
\widetilde X_0 \to \widetilde \Sigma$ the projections. 

\smallskip \noindent We first claim that $\tilde r_0$ is an elliptic fibration. Indeed :
$$K_{\widetilde X_0} =
g^\ast K_{X_0} +  B, \eqno (*)$$ where
 where $B$ is
the exceptional set of the generically finite map $g$ (recally that there is no ramification in 
codimension 1). 
Let $\tilde F$ be a general fiber of $\tilde r_0$ and $F$ a general fiber of
$r_0.$ Since $K_{X_0} \cdot F = 0$ we conclude by (*) that $K_{\tilde X_0}Ê\cdot \tilde F = 0,$
hence $\tilde F$ is elliptic. 

\noi Now blow up $\tilde X_0$ in order to make $p_2$ holomorphic, if necessary ($p_1$ is anyway
holomorphic, since $\tilde E$ is elliptic).

\noi We next claim the existence of a map $\pi : \widetilde S \to
\widetilde E$ such that
$\pi \circ \tilde r_0 = p_1 : \widetilde X_0 \to \widetilde E$.
This comes down to show that, if
$\widetilde F_0$ is a generic fiber of $\tilde r_0$ then  $ p_1(\widetilde F_0) \not=
\widetilde E$. Otherwise, every such $\widetilde F_0$ is a finite Žtale
cover of
$\widetilde E$, and thus every fiber $F_0$  of $r_0$ is isogeneous to
$\widetilde E$. This
would imply, that for $P
\subset S$ a rational
curve, the preimage $r^{-1}_0(P)$ is a birationally ruled elliptic surface, and
sothe rationality of $S$ yields
$\kappa(X_0) = -\infty$, contrary to our hypothesis.
Now observe that $\sigma: \widetilde S \to S$ is finite with finite ramification locus,
i.e. \'etale in codimension 2. This implies easily that $\widetilde S$ is rational, too,
e.g. by observing the following 
commutative diagram :
$$\diagram{
\widetilde S &\hfl{}{} &S \cr
\vfl{\pi}{} & &\vfl{}{\pi'} \cr
\widetilde E &\hfl{u}{} &E \cr}$$

with $u$ Žtale and $E$ elliptic. Or by observing that there are a lot of rational curves
on $\widetilde S$ by lifting rational curves from $S.$ In total we obtain a 
contradiction to the existence of $\pi,$ which yields a 1-form on $\tilde S.$ 

\bn {\bf (8.7)} $\underline{a(X) = 1}$. We again consider the algebraic
reduction $r : X \to
C$ of $X$, where $C$ is now a curve. Let $F$ be a general fiber of $r.$ When
$F$ is a torus, we let $r_0: X_0 \to C$ again be the untwisted model of $r.$ Finally we let 
$ a(r) := a(F)$.

\noi We consider successively the cases $a(r) = 2, 1, 0$ in 8.8, 8.9 and
8.10 respectively.

\bn {\bf (8.8)} $\underline{a(X) = 1 \, ; a(r) = 2}$. Again we have two
cases: $F$
is either a torus or ${\Bbb P}({\cal O}_E \oplus {\cal L})$ as in (7.3.ii). In
the last case, we have $\kappa(X) = -\infty$ and $X$ is uniruled, as
claimed. So we have
only to concentrate on the first case.

\bn {\bf (8.8.1) Claim :} $\underline{\chi({\cal O}_X) = 0 \, ; q(X) > 0}$.
\smallskip \noi {\bf Proof.} Let $D \subset X_0$ be the canonical section of $r_0.$
Let $r'_0: X'_0 \la C$ be a relative minimal model, so $K_{X'_0}$ is $r'_0-$nef and 
moreover $mK_{X'_0}$ is $r'_0-$generated [KMM87,6-1-13], i.e.
$$ r'^*_0 {r'_0}_*(mK_{X'_0}) \la mK_{X'_0}Ê$$
is surjective.
Since ${r'_0}_*(mK_{X'_0})$ is a line bundle, we conclude 
$$ r'^*_0 {r'_0}_*(mK_{X'_0}) \simeq mK_{X'_0},$$
so that $mK_{X'_0}Ê\simeq r'^*_0(L')$ for some line bundle $L'$ on $C.$ 
Since $r_0$ has a section $D,$ it follows that $r'_0$ cannot have multiple fibers, all
fibers have to be reduced. Then we conclude immediately that $K_{X'_0} \simeq {r'_0}^*(L)$
for some line bundle $L$ on $C.$ So $X'_0$ is Gorenstein and by Riemann-Roch [Fl87],
we obtain
$$ \chi(X'_0,\cO_{X'_0}) = - {1 \over {24}}K_{X'_0}Ê\cdot c_2(X'_0).$$
Let $G$ be the general fiber of $r'_0.$ Then $K_{X'_0} \equiv kG$ for some integer $k.$
Therefore $c_2(X'_0) \cdot G = c_2(G) = 0,$ since $G$ is a torus. Thus $\chi(X'_0) = 0,$
so $$\chi(X,\cO_X) = \chi(X_0,\cO_{X_0}) = 0.$$
Since $h^2(\cO_X) > 0$ and since $\kappa (X) = 0,$ we finally conclude $q(X) \geq 1,$ proving
the claim (8.8.1).

\medskip

\noi {\bf (8.8.2) Subcase :} $\underline{\kappa(X) = -\infty}$. Consider the
Albanese map
$\alpha : X \to \hbox{Alb}(X)$ of $X$ and conclude from either $C_{3,2}$ or
$C_{3,1}$ that
its fibers have $\kappa = -\infty$, so that $X$ is uniruled.

\medskip

\noi {\bf (8.8.3) Subcase :} $\underline{\kappa(X) = 0}$. In this case we
may have $q(X)
= 1, 2, 3$. We subdivide again in three subcases successively: 

\medskip
\noi {\bf (8.8.3.1)} $\underline{q(X) = 1}$. 
\sn  Notice that $r$ is locally Moishezon. On the other hand, $X$ is K\"ahler and
therefore $r$ is locally projective [10.1].   
By [Ka88], $r$ admits locally a relative minimal model via relative contractions and flips. 
Since the general fiber of $r$ is a torus, all local contractions and flips glue to global
ones whence $r$ admits globally a relative minimal model $r' : X' \la C.$
We claim that $mK_{X'} \simeq \cO_{X'}.$ Write $kK_{X'}Ê\simeq \cO_{X'}(E)$ for some positive
$k.$ Let $F'$ be a general fiber of $r'.$ By adjunction we have $K_{X'}Ê\vert F' = \cO.$
Thus $E \cap F' = \emptyset$ and therefore ${\rm dim}r'(E) = 0.$ Since $E$ is relatively nef,
$E \equiv rF'$ for some positive rational $r.$ This means that a multiple of $E$ consists only
of fibers of $r',$ hence $mK_{X'}Ê\simeq r^{'*}(L).$ Since $\kappa (X') = 0,$ the line bundle
$L$ must be trivial. So $mK_{X'}Ê\simeq \cO.$ Take a canonical cover of $X'$Ê[KMM87]; this is
a covering $\hat X \la X'$ which is \'etale outside the singularities of $X',$ hence \'etale 
in codimension 2. Take a non-zero 2-form $\eta$ on $\tilde X$ (= image of a 2-form on a 
desingularisation). Then by [Pe94] there is a 1-form $\omega$ on $X'$ such that $\omega 
\wedge \eta$ is a non-zero, hence nowhere vanishing 3-form on $X',$ i.e. a section in $K_{X'}.$ 
In particular $\omega$ has no zeroes on the regular part $\hat X_0$. Hence $\hat r: \hat X_0
\la C$ is a submersion. Now consider a singular fiber $D$ of $\hat r.$ Since $D$ has only
finitely many singularities, it is normal. Since $\omega_D \simeq \cO_D;$ we have 
$H^1(\cO_D) = 0$ by [Um81]. Since $h^1(\cO_{\hat F}) = 2$ for the general fiber $\hat F$
of $\hat F,$ this contradicts the semi-continuity theorem, $\hat r$ being flat. Thus 
$\hat r$ is a submersion, and $\hat X$ is smooth. Now it is classical (e.g. [Be83]) that
$\hat X$ is a product, possibly after another finite \'etale cover. 

\medskip
\noi {\bf (8.8.3.2)} $\underline{q(X) = 2}$.
\sn In (8.8.3.1) we used the assumption $q(X) = 1$ only in order to conclude that 
$C$ is elliptic (since $\kappa (X) = 0$, we cannot have $g(C) \geq 2$). So suppose $C$
rational. Then the Albanese map $\alpha : X \to A$ must be surjective and $a(A) = 0.$ 
On the other hand the fibers of $r$ are projective and map onto $A.$ This is a contradiction.
Thus $C$ is elliptic and we can apply the same arguments as in (8.8.3.1).
\medskip \noi {\bf (8.8.3.3)} $\underline{q(X) = 3}$. Because $\kappa(X) = 0$,
$\alpha : X \to
\hbox{Alb}(X)$ is bimeromorphic.

\medskip

\noi {\bf (8.9)} $\underline{a(X) = a(r) = 1}$. This case is very similar
to the next one
treated in 8.10 below. We give however in case $\kappa(X) = -\infty$  an
alternative short
proof :

\noi Let $g : X \to Y$, $h : Y \to C$ with $h \circ g= r$ the relative
algebraic reduction of
$r$. Then $Y$ is a surface with $a(Y) = a(X) = 1$ and with algebraic reduction
$h$. In
particular, $\kappa(Y) \ge 0$ and we conclude from $C_{3,2}$ that $X$ is
uniruled.

\bn {\bf (8.10)} $\underline{a(X) = 1 \, ; a(r) = 0}$. Then $F$, the
generic fiber of $r$, is
isomorphic to a fixed surface which is either a torus or a K3-surface (use 7.5).
Assume first that
$F$ is a torus and consider the untwisted model $r_0 : X_0 \to C.$ By (6.7) we may assume that 
$r_0 : X_0 =  F
\times \widetilde C)/G \to C = \widetilde C/G$ for some finite group $G$ acting
diagonally on $F \times \widetilde C$.

\noi We therefore have two maps
$$\diagram{
X_0 &\hfl{r_0}{} &C \cr
\vfl{g}{} & & \cr
Z := (F/G) & & \cr}$$
Notice that $h^0(\Omega^2_Z) \ne 0$ since $a(Z) = a(F) = 0.$ 
First suppose $\kappa(X) = -\infty.$ By (6.10), the canonical map $X_0 \to C \times Z$ is
bimeromorphic. Hence the restriction $F \to Z$ is bimeromorphic, and thus $q(X_0) > 0.$ 
Hence $q(X) > 0$ by (6.7), and $X$ is uniruled via the Albanese map.

\noi Assume now that $\kappa(X) = 0$ ; then $q(C) \le 1$ by $C_{3,1}$.
Assume first that
$C = {\Bbb P}_1$. Argiunf as before, we have $q(X) = 2.$ Then the
Albanese map
$\alpha : X \to \hbox{Alb}(X)$ is a connected surjective elliptic
fibration.

\noi As in (8.8.3.2) we conclude that there is a finite cover $X'Ê\to X,$ \'etale in codimension 1,
such that $X'$ is bimeromorphically a torus. But this is only possible if $C$ is elliptic.
So this situation
does not occur.

\noi Assume now that $\kappa(X) = 0$ and that $C$ is elliptic.
The fibers of $g$ are elliptic curves since $\kappa(Z) = \kappa(X_0) =
0$. The fibers of
$g$ are thus mapped to $C$ in an Žtale way by $r_0$. Via the finite \'etale base
change
$\beta : \overline C \to C$, the map $\overline r_0 \times \overline g :
\overline X_0 \to
Z \times \overline C$ is bimeromorphic; here $\overline X_0 = X_0 \dis
\mathop
{\times}_{C} \overline C$. But then $q(\overline X_0) = 3$ and $\overline
X_0$ is
bimeromorphic to a torus. Considering $\overline X := X \dis \mathop
{\times}_{C}
\overline C$ (which is Žtale over $X$), we get $q(\overline X) = 3$ and the
conclusion.

\medskip \noindent We still have to treat the case where $F$ is a K3-surface (with $a(X)
= 1$, $a(r) = 0$). By (7.6) we have a meromorphic map to a K3-surface. Therefore $X$ is
uniruled by $C_{3,2}$ if $\kappa (X) =  -\infty. $ For the case $\kappa (X) = 0,$
we proceed as in (8.8.2.2).

\medskip \noindent This concludes the proof of 8.1.

\bn \bn \bn  {\medium 9. Classification of non-algebraic threefolds with  $\kappa =
-\infty$.}
\bn {\bf 9.1 Theorem} {\it Let $X$ be a compact connected non-algebraic
K\"ahler
threefold with $\kappa(X) = -\infty$. Then $X$ is bimeromorphic to exactly
one of the
following :\hfill\break
{\rm (i)} $\underline{a(X) = 0} :$ $X$ is simple but not Kummer.
\smallskip
\noi{\rm (ii)} $\underline{a(X) = 0} :$ $X$ is not simple. Then it has a
unique map $\rho : X
\to S$ to a surface $S$ with generic fiber ${\Bbb P}_1$. Moreover $S$ is bimeromorphically a torus or
a K3-surface
with $a(S) = 0$.
\smallskip
\noi {\rm (iii)} $\underline{a(X) = 1 \, ; a(r) = 2} :$ The generic fiber
of $r : X \to C$ (the
algebraic reduction of $X$) is
\smallskip
\noi ${\Bbb P}({\cal O}_E \oplus {\cal L})$ as in (7.3.ii). Let
 $f : X \to
S$, $g : S \to C$ with $g \circ f = r$ be the relative Albanese map of $r.$ Then $a(S) = 1$ and $g$ is the
algebraic reduction of
$S$. The generic fiber of $f$ is ${\Bbb P}_1$.
\smallskip
\noi {\rm (iv)} $\underline{a(X) = 1 \, ; a(r) = 0} :$ $X$  is bimeromorphically $ {\Bbb P}_1
\times F$ with $F$ either
a torus or a K3-surface with $a(F) = 0$.
\smallskip
\noi {\rm (v)} $\underline{a(X) = 2} :$ There exists a (unique) map $\sigma
: X \to \Sigma$
with generic fiber ${\Bbb P}_1$ and
$a(\Sigma) = 1$. Moreover there is a ruling $\rho : S \to
C$ and an algebraic reduction $\tau : \Sigma \to C$ such that the product
map $\sigma
\times r : X \to \Sigma \dis \mathop {\times}_{C} S$ is onto.}

\noi {\bf Proof.} This has been established during the proof of 8.1; for (iv) use also (7.6) and
(6.10).

\bn \bn {\medium 10. A projectivity criterion for Moishezon morphisms}

\bn In this last section we prove the following result which was already used in section 8.
Recall that a proper morphism is Moishezon if it is bimeromorphically equivalent to a
projective morphism.

\bn {\bf 10.1 Theorem}  {\it Let $f: X \la S$ be a surjective Moishezon morphism between complex manifolds

$X$ and $S.$ If $X$ is K\"ahler, 
$f$ is locally (over $S$) projective. }

\bn If $f: X \la S$ is a proper morphism of complex manifolds, an $f-$ curve is by
definition an irreducible compact curve which is contained in some fiber of $f.$ We
let $N_1(X/S) \subset H_2(X,{\bf Q})$ denote the subspace generated by the classes of
$f-$curves and $N_1^*(X/S)$ its dual (over ${\bf Q}).$ In this terminology Theorem 10.1 will be a consequence of the following

\bn {\bf 10.2 Theorem} {\it Let $f: X \la S$ be a surjective Moishezon morphism of complex manifolds.
Let $s \in S$ and $\Lambda \in N_1^*(X/S).$ Then - possibly after shrinking $S$ around $s$ -
there exists a line bundle $L \in {\rm Pic}(X)$ such that a positive rational number $m$ such that 
$m\Lambda = L,$ i.e.
$$ m\Lambda (C) = L \cdot C$$
for all $f-$curves $C.$}

\bn In the absolute case ($S$ a point) both statements are due to Moishezon [Ms67].
Usually $N_1(X /S)$ is defined in another way,  by taking numerical equivalence using holomorphic
line bundles, see e.g. [KMM87]. It follows from (10.2) that both notions coincide. 

\bn {\bf(10.3)} Here we show how to derive (10.1) from (10.2).
\sn Fix $s \in S$ and a K\"ahler form on $X.$ We define $\Lambda \in N_1^*(X/S)$ by
$$ \Lambda (C) = \int_C \omega $$
for all irreducible $f-$curves $C \subset X,$ hence all $C \in N_1(X/S).$ Then for given 
$\epsilon > 0$ we can find
$\Lambda_{\epsilon} \in N_1^*(X/S)$ such that 
$$ \vert \Lambda (C) - \Lambda_{\epsilon}(C) \vert \leq \epsilon \Lambda (C).$$ 
Applying Theorem 10.2 to $\Lambda_{\epsilon}$ we find $m > 0$ and a divisor $L_{\epsilon}$
such that $m \Lambda_{\epsilon} = L_{\epsilon},$ hence 
$$ \vert m \Lambda (C) - L_{\epsilon}Ê\cdot C \vert \leq m \epsilon \Lambda (C) \eqno (*)$$
for $C \in N_1(X/S).$ 
We shall show that $L_{\epsilon}$ is $f-$ample, provided $\epsilon > 0$ is sufficiently
small. From (*) we deduce that $L_{\epsilon}$ is $f-$nef if $\epsilon \leq 1$ and that 
$L_{\epsilon}Ê\vert X_t$ is ample for all $t \in S$ by Kleiman's criterion. Now the
claim follows from

\bn {\bf 10.4 Lemma} \ ÊÊ{\it Let $X$ and $S$ be reduced complex spaces, $f: X \la S$
a proper surjective map. Let $L$ be a line bundle on $X$ such that $L \vert X_s$ is
ample for every $s.$ Then $f$ is locally (over $S$) projective. }

\bn {\bf Proof.} Since the claim is local over $S,$ we may assume $S$ Stein. We proceed
by induction over $d = \dim X.$ If $d = 0,$ the claim is obvious; so suppose $d > 0.$ By
substituting $L$ by a multiple we may assume $f_*(L) \ne 0,$ hence, $S$ being Stein,
$H^0(L) \ne 0.$  
We check the  relative ampleness of $L$ (locally over $S$) by the relative Serre vanishing
$$ R^pf_*(\F \otimes L^{m}) = 0$$
for $p > 0,$ every coherent sheaf $\F$ and $m \geq m_0(\F).$ 
So fix a coherent sheaf $\F$ on $X$ and denote $S_m$ the support of the coherent sheaf
$R^pf_*(\F \otimes L^m)$ for a fixed $p.$ We claim that 
$$ S_{m+1}Ê\subset S_m \eqno (*)$$
for sufficiently large $m.$ In fact, take $D \subset \vert L^m \vert.$ Then we have a
sequence
$$ 0 \la L^m \la L^{m+1}Ê\la L^{m+1}Ê\vert D \la 0.$$
Now, taking direct images and applying induction, claim (*) is clear. Fix $s_0 \in S$. 
It follows that for all large $m, $ the point $s_0$ is not contained in $S_m.$ Moreover,
since $S_{m+1}Ê\subset S_m,$ there is a neighborhood $U$ of $s_0$ disjoint from all
these $S_m.$ Substituting $S$ by $U$ we have the relative vanishing we were looking for. 

\bn {\bf Proof of 10.2.}
\sn (1) First we notice that it is sufficient to prove the theorem for $f$ projective; of course
we may always assume $S$ Stein. To reduce to the projective case, let $f: X \la S$ be
a Moishezon morphism. Let $g: \hat X \la X$ be a sequence of blow-ups with compact smooth
centers (in fibers of $f$) such that the induced morphism $\hat f: \hat X \la S$ is projective.
Let $$\hat \Lambda = \Lambda \circ g_* : N_1(\hat X /S) \la {\bf Q}.$$ Then by assumption
there exists $\hat D \in {\rm Pic}(\hat X)$ and $m > 0$ such that $m \hat \Lambda = \hat D.$
Let $D = g_*(\hat D).$ Then it is immediately checked that $m \Lambda = D.$ 
\sn So from now on we shall assume that $ f$ is projective and $S$ is Stein. Since 
our claim is local we fix a point $s_0 \in S$ and may always shrink around $s_0.$ 

\sn (2) We assume here that $R^2f_*(\O_X)$ is torsion free. This is e.g. guaranteed if there is a normal
crossings divisor $H \subset S$ such that $f$ is smooth over $S \setminus H,$ see 
[Ko86],[Mw87],[Na86].
Fix $\Lambda' \in N_1^*(X/S).$ Then there is an  
element $\Lambda \in H^2(X,{\bf Q})$ such that $\Lambda \cdot C = \Lambda'(C)$ for all
curves $C$ contracted by $f.$ After multiplying $\Lambda'$ by a suitable positive integer,
we may assume that $\Lambda \in H^2(X,{\bf Z}).$ 
We consider the following commutative diagram, all maps being canonically defined
$$ \matrix { H^2(X,{\bf Z}) & \buildrel {\psi}Ê\over \laÊ& H^2(X,\O_X) Ê\cr
u \downarrow &  & \downarrow v   \cr
H^0(S,R^2f_*({\bf Z})) & \buildrel {\psi'}Ê\over \la & H^0(S,R^2f_*(\O_X)) \cr }.$$ 
Since  $S$ is Stein, $v$ is an isomorphism. 
We construct a divisor $D'  \in {\rm Pic}(X)$ such that over a  Zariski open dense set $ S' \subset S^*$  we have 
$$v \circ \psi (\Lambda - c_1(D')) = 0.  \eqno (*)$$ 
Given (*), it follows $v\psi(\Lambda - c_1(D')) = 0$ over $S$, since $R^2f_*(\O_X) $ is torsion
free. Therefore $\psi (\Lambda) =\psi(\Lambda - c_1(D')) = 0 $ and hence $\Lambda = c_1(L)$ on $X$ for some 
$L \in {\rm Pic}(X).$  
To prove (*) notice first that on by [Ms67] there is $D_s \in {\rm Pic}(X_s)$ such that 
$$c_1(D_s) = \Lambda$$
for all $s \in S^*.$ Fix $H \in {\rm Pic}(X)$ $f-$ample. Then for every $s \in S^*$ we find positive
integers $N_s$ and $\mu_s$ such that 
$$ \mu(N_sH - D_s)$$
is very ample on $X_s.$ Moreover we can choose $\mu$ and $N$ independent on $s$ for
$s$ in a non-empty Zariski-open subset $S' \subset  S^*.$ Let $m = {\dim}  X - {\dim} S.$ 
Choose an irreducible component $G \subset {\cal C}_{m-1}(X/S)$ of the relative cycle space
containing $m-1-$cycles of $X$ contained in some fiber of $f$ and such that for generic
$Z \in G$ the Cartier divisor $Z \subset X_s$ satisfies 
$$ c_1(Z) = c_1(\mu(NH - D_s)).$$ 
Let $f_* : G \la S$ be the map associating to a cycle $Z$ the point $s$ with $Z \subset f^{-1}(s).$
Then we can choose $G$ in such a way that $f_*$ is surjective and proper and moreover that there
is a closed analytic set $\tilde S \subset G$ satisfying
{\item {(a)} $f_*(\tilde S) = S$ and
\item {(b)}Ê$f_* \vert \tilde S \la S$ is generically finite.}
Possibly we have to shrink $S$ around $s_0.$ 
Let $\tilde D$ be the graph of the family of $f-$cycles parametrised by $\tilde S $
and let $D$ be its image in $X.$ Then $D$ is a divisor in $X$ such that 
$$ D \cdot C_s = \mu(NH - D_s) \cdot C_s$$
for all $s \in S'$ and all irreducible curves $C_s \subset X_s.$ In other words:
$$ \Lambda \cdot C_s = D' \cdot C_s.$$
with $D' = \mu(NH _ D).$ Therefore $u(\Lambda - c_1(D')) = 0$ over $S'$ and our commutative diagram
above gives claim (*).
\sn (3) In general the arguments in (2) show that we can still construct $D'$ such that $D' \cdot C_s = \Lambda \cdot C_s$ for all
$s \in S'.$ Let 
$$ \sigma = \psi'u(\Lambda - c_1(D')) = \psi'u(\Lambda) \in H^0(S,R^2f_*(\O_X)).$$
Let $T = {\rm Supp}(\sigma).$ If $T = \emptyset,$ we conclude by the old arguments. So suppose
$T \ne \emptyset;$ this case has to be ruled out.  
\sn (4) Let $A = f^{-1}(T)$ and $U = X \setminus A.$
We have $\psi' u(\Lambda) = 0$ on $S \setminus T.$ Therefore $v\psi(\Lambda) = 0$ on $S \setminus T.$
It follows from the functoriality of the Leray spectral sequence that actually
$\psi(\Lambda) = 0$ on $X \setminus A.$ So the restriction map 
$$ H^2(X,{\bf Z}) \la H^2(U,{\bf Z})$$
maps $\Lambda$ to $0$ and thus we obtain 
$$\Lambda \in {\rm Im}(H^2(X,U) \la H^2(X)).$$ 
By Lemma 10.5 below we see that, computing over $\Bbb Q,$ 
$$ H^2(X,U) \simeq H_{2n-2}(A) \simeq \bigoplus {\Bbb Q}Ê[A_i], $$
where the $A_i$ vary over the irreducible components of $A$ of codimension 1 
in $X.$ It follows $r\Lambda \in {\rm Pic}(X)$ fo a suitable positive integer,
hence $\lambda \in {\rm Pic}(X)$ (showing also that $T$
is empty). This ends the proof of (10.2).

\bn {\bf 10.5 Lemma}  {\it Let $f: X \la S$ be a proper surjective map of complex manifolds with connected
fibers. Let $n = {\dim }X.$ Let $T \subset S$ be an analytic set, $A = f^{-1}(T)$ and $U = X \setminus A.$ Then
$$ H^2(X,U) \simeq H_{2n-2}(A) \simeq \bigoplus {\Bbb Q}Ê[A_i], $$
where the $A_i$ vary over the irreducible components of $A$ of codimension 1 in $X.$}

\bn {\bf Proof.}
Let $B \subset A$ be the singular locus of $A,$ $\pi: \hat A \la A$ a 
desingularisation and $\hat B = \pi^{-1}(B).$ First notice that
$$H^{2n-2}_c(A,B) \simeq H^{2n-2}_c(\hat A,\hat B). \eqno (1)$$
In fact, this follows from
$$ H^{2n-2}_c(A,B) \simeq H_0(A\setminus B)$$
(and analogously for $(\hat A,\hat B)$) which is a consequence of [Sp66,thm19,p.297] via one point
compactification. 
Since $H^{i}_c(\hat B) = H^{i}_c(B) = 0 $ for $i \geq 2n-1,$ we conclude from (1):
$$ H^{2n-2}_c(\hat A,) \simeq H^{2n-2}_c(A). \eqno (2)$$ 
Since
$$ H^{2n-2}_c(\hat A) = \bigoplus \Bbb Q [A_i], $$
the sum taken over all $(n-1)-$dimensional components of $\hat A,$ (2) gives
$$ H^{2n-2}_c(A) \simeq \bigoplus \Bbb Q [A_i]. \eqno (3)$$
Finally, we have
$$ H^{2n-2}_c(A) \simeq H_2(X,U)$$
by [Sp66,thm10,p.342], hence the universal coefficient theorem yields
$$ H^{2n-2}_c(A) \simeq H^2(X,U),$$
proving in combination with (3) our claim .

\vfill \eject \noindent {\medium References}

\bn 

\item {[Be83]}ÊBeauville,A.: Vari\'et\'es k\"ahl\'eriennes dont la premiere classe de Chern
est nulle. J. Diff. \ \   Geom. 18, 755-782 (1983)

\item {[Bl58]} Blanchard, A.: Sur les vari\'et\'es analytique complexes. Ann. Sci. Ec. Norm. Sup. 73,
157-202 (1958)

\item {[BPV84]} {\obeylines Barth,W.;Peters,C;van den Ven, A.: Compact complex surfaces. Erg. d. Math.,
3.Folge, Band 4. Springer 1984}

\item {[Ca81]} Campana,F.: Cor\'eduction alg\'ebrique d'un espace analytique faiblement
k\"ahl\'erien compact. Inv. math. 63, 187 - 223 (1981)

\item {[Ca85]}Ê{\obeylines Campana,F.: R\'eduction d'Albanese d'un morphisme propre et faiblement K\"ahl\'erien. 
Comp. Math.54, 373-416 (1985)}

{\obeylines \item {[Ca92]} Campana,F.: Connexit\'e rationnelle des vari\'et\'es de Fano. Ann. scient.
Ec. Norm. Sup. 25, 539-545 (1992)

\item{[Fl87]} Fletcher,A.R.: Contributions to Riemann-Roch on projective 3-folds with only 
canonical singularities. Proc. Symp. Pure Math. 46, 221-231 (1987)}

\item {[Fu83]} Fujiki,A.: On the structure of compact complex manifolds in class ${\cal C}.$
Adv. Stud. Pure Math. 1, 231 - 302 (1983)

\item {[Fj78]}ÊFujita,T.: K\"ahler fiber spaces over curves. J. Math. Soc. Japan 30, 779 - 794
(1978)

\item {[Ka81]} Kawamata,Y.: Characterisation of abelian varieties. Comp. math. 43, 253 - 276
(1981)

\item {[Ka85]} Kawamata,Y.: Minimal models and the Kodaira dimension of algebraic
fiber spaces. J. reine u. angew. Math. 363, 1 - 46 (1985)

\item {[Ka88]} {\obeylines Kawamata,Y.: The crepant blowing-ups of 3-dimensional canonical 
singularities and its application to degeneration of surfaces.
Ann. Math. 127,93-163 (1988)}

\item {[KMM87]} Kawamata,Y.;Matsuda,K.;Matsuki,K. : Introduction to the minimal model problem.
Adv. Stud. Pure Mat. 10, 283-360 (1987)

\item {[Ko86]} Koll\'ar,J.: Higher direct images of dualizing sheaves II. Ann. Math. 124, 171-202 (1986)

\item {[KoMiMo92]} Koll\'ar,J.;Miyaoka,Y.;Mori,S.: Rational connectedness and boundedness for
Fano manifolds. J. Diff. Geom. 36, 765 - 769 (1992)

\item {[Li78]} Lieberman,D. : Compactness of the Chow scheme. Lect. Notes in Math. 670, 140-186 (1978)

\item {[Mi87]} Miyaoka,Y.: Deformation of a morphism along a foliation and applications. Proc. Symp.
Pure Math. 46, 245-268 (1987)

\item {[Mi88]} Miyaoka,Y.: On the Kodaira dimension of minimal threefolds. Math. Ann. 281, 325-332
(1988)

\item {[Mo82]} Mori,S.: Threefolds whose canonical bundles are not numerically effective. 
Ann. Math. 116, 133 - 176 (1982) 

 \item {[Mo88]} Mori,S.: Flip theorem and the existence of minimal models for 3-folds.
J. Amer. Math.   Soc. 1, 117-253 (1988)

\item {[Ms67]}ÊMoishezon,B.G.: On $n-$dimensional compact varieties with $n$ algebraically
independent meromorphic functions. Amer. Math. Soc. Transl. 63, 51-177 (1967) 

\item {[Mw87]} Moriwaki,A.: Torsion freeness of higher direct images of canonical bundles,
Math. Ann. 276, 385-398 (1987)

\item {[Na86]} Nakayama,N.: Hodge filtration and the higher direct images of dualising
sheaves. Inv. math. 86, 217-221 (1986)

\item {[Na88]} Nakayama,N.: On Weierstrass models. Alg.Geom. and Comm. Algebra; vol. in
honour of Nagata, vol. 2, 405 - 431. Kinokuniya, Tokyo 1988

\item {[Na95]} Nakayama,N.: Local structure of an elliptic fibration. Preprint 1995

\item {[Pe94]} Peternell,T.: Minimal varieties with trivial canonical class,I. Math. Z.
217, 377 - 407 (1994)

\item {[Pe96]} {\obeylines Peternell,T.: Towards a Mori theory on compact K\"ahler threefolds, II.
Preprint 1996, to appear in Math. Ann. }

\item {[Sp66]}ÊSpanier,E.: Algebraic topology. McGraw-Hill 1966

\item {[Ue75]} Ueno,K.: Classification theory of algebraic varieties and compact complex
spaces. Lecture Notes in Math. 439. Springer  1975

\item {[Ue87]} Ueno,K.: On compact analytic threefolds with non-trivial Albanese torus.
Math. Ann. 278, 41 - 70 (1987)

\item {[Um81]} Umezu,Y.: On normal projective surfaces with trivial dualising sheaf. Tokyo J.
Math. 4, 343-354 (1981)

\item {[Vi80]} Viehweg,E.: Klassifikationstheorie algebraischer Variet\"aten der Dimension 3.
Comp. math. 41, 361-400 (1980)

\vskip .5cm
\vskip .5cm
$$\matrix { {\rm  Frededric \  Campana} & & & & &  {\rm Thomas \  Peternell} \cr
{\rm Departement \  de  \ mathematiques} & & & & & {\rm Mathematisches \  Institut} \cr
{\rm Universite \  de \ Nancy} & & & & & {\rm Universitaet \  Bayreuth} \cr
{\rm BP 239} & & & & & & \cr
{\rm 54506 \  Vandoeuvre \  les \  Nancy} & & & & & {\rm 95440 \  Bayreuth} \cr
{\rm France} & & & & & {\rm Germany} \cr}Ê$$

\end